\documentclass{amsart}

\usepackage{amsfonts}  
\usepackage{amssymb}
\usepackage{enumitem}

\usepackage{moreenum}
\usepackage{amsmath}
\usepackage[T1]{fontenc}
\usepackage{graphicx}
\usepackage{lmodern}
\usepackage[unicode]{hyperref}

\newtheorem{ttt}{Theorem}[section]
\newtheorem{llll}[ttt]{Lemma}
\newtheorem{ccc}[ttt]{Claim}
\newtheorem{eee}[ttt]{Example}
\newtheorem{fff}[ttt]{Fact}
\newtheorem{rrr}[ttt]{Remark}
\newtheorem{sss}[ttt]{Statement}
\newtheorem{ddd}[ttt]{Definition}
\newtheorem{qqq}[ttt]{Question}
\newtheorem{cccc}[ttt]{Corollary}
\newtheorem{nnn}[ttt]{Notation}
\newtheorem{ppp}[ttt]{Problem}
\newtheorem{pppp}[ttt]{Proposition}
\newtheorem{ccccc}[ttt]{Conjecture}

\newcommand{\beq}{\begin{equation} }

\newcommand{\bt}{\begin{ttt}}
\newcommand{\bl}{\begin{llll}}
\newcommand{\bc}{\begin{ccc}}
\newcommand{\bex}{\begin{eee}}
\newcommand{\bfa}{\begin{fff}}
\newcommand{\br}{\begin{rrr}\upshape}
\newcommand{\bst}{\begin{sss}}
\newcommand{\bd}{\begin{ddd}\upshape}
\newcommand{\bdd}{\begin{ddd}\upshape}

\newcommand{\bq}{\begin{qqq}}
\newcommand{\bnn}{\begin{nnn}}
\newcommand{\bpr}{\begin{ppp}}
\newcommand{\bprop}{\begin{pppp}}
\newcommand{\bcor}{\begin{cccc}}
\newcommand{\bcon}{\begin{ccccc}}

\newcommand{\eeq}{\end{equation}}
\newcommand{\et}{\end{ttt}}
\newcommand{\el}{\end{llll}}
\newcommand{\ec}{\end{ccc}}
\newcommand{\eex}{\end{eee}}
\newcommand{\efa}{\end{fff}}
\newcommand{\er}{\end{rrr}}
\newcommand{\est}{\end{sss}}
\newcommand{\ed}{\end{ddd}}
\newcommand{\eq}{\end{qqq}}
\newcommand{\ecor}{\end{cccc}}
\newcommand{\econ}{\end{ccccc}}
\newcommand{\enn}{\end{nnn}}
\newcommand{\epr}{\end{ppp}}
\newcommand{\eprop}{\end{pppp}}

\newcommand{\bp}{\noindent\textbf{Proof. }}
\newcommand{\ep}{\hspace{\stretch{1}}$\square$\medskip}

\newcommand{\sm}{\setminus}

\newcommand{\beeq}{\begin{equation}}
\newcommand{\eeeq}{\end{equation}}

\newcommand{\cf}{\mathrm{cf}} 

\newcommand{\dom}{{\rm dom}}
\newcommand{\ran}{{\rm ran}}

\newcommand{\bref}[1]{\boxplus_{\ref{#1}}}

\newcommand{\up}{\upharpoonright}

\DeclareMathOperator{\supp}{supp}

\DeclareMathOperator{\htt}{ht}

\DeclareMathOperator{\otp}{otp}

\DeclareMathOperator{\id}{id}
\DeclareMathOperator{\trcl}{trcl}
\DeclareMathOperator{\ddf}{def}

\DeclareMathOperator{\Sp}{Sp}

\newcommand{\tieconcat}{%
	\mathbin{\mathpalette\dotieconcat\relax}%
}
\newcommand{\dotieconcat}[2]{
	\text{\raisebox{.8ex}{$\smallfrown$}}%
}

\numberwithin{equation}{section}

\def\name#1{\mathchoice%
	{\setbox0=\hbox{$\displaystyle #1$}
		\setbox1=\vtop{\ialign{##\crcr
				$\hfil{\displaystyle #1}\hfil$\crcr\noalign{\kern2pt%
					\nointerlineskip}$\hfil\mathord{\displaystyle \sim}%
				\hfil$\crcr\noalign{\kern2pt\nointerlineskip}}}%
		\vphantom{\copy1}%
		\setbox2=\hbox{$\displaystyle \sim$}%
		\wd1=\wd0\dp1=0cm\ifdim\wd2>\wd1 \wd1=\wd2\else\relax\fi
		\ht1=\ht0\relax
		\box1}%
	{\setbox0=\hbox{$\textstyle #1$}
		\setbox1=\vtop{\ialign{##\crcr
				$\hfil{\textstyle #1}\hfil$\crcr\noalign{\kern1.2pt%
					\nointerlineskip}$\hfil\mathord{\textstyle \sim}%
				\hfil$\crcr\noalign{\kern1.5pt\nointerlineskip}}}%
		\vphantom{\copy1}%
		\setbox2=\hbox{$\textstyle \sim$}%
		\wd1=\wd0\dp1=0cm\ifdim\wd2>\wd1 \wd1=\wd2\else\relax\fi
		\ht1=\ht0\relax
		\box1}{\setbox0=\hbox{$\scriptstyle #1$}
		\setbox1=\vtop{\ialign{##\crcr
				$\hfil{\scriptstyle #1}\hfil$\crcr\noalign{\kern1pt%
					\nointerlineskip}$\hfil\mathord{\scriptstyle \sim}%
				\hfil$\crcr\noalign{\kern2.1pt\nointerlineskip}}}%
		\setbox2=\hbox{$\scriptstyle \sim$}%
		\vphantom{\copy1}%
		\wd1=\wd0\dp1=0cm\ifdim\wd2>\wd1 \wd1=\wd2\else\relax\fi
		\ht1=\ht0\relax
		\box1}{\setbox0=\hbox{$\scriptscriptstyle #1$}
		\vtop{\ialign{##\crcr
				$\hfil{\scriptscriptstyle #1}\hfil$\crcr\noalign{\kern1pt%
					\nointerlineskip}$\hfil\mathord{\scriptscriptstyle \sim}%
				\hfil$\crcr\noalign{\kern1.5pt\nointerlineskip}}}%
}}

\begin{document}

\author[M\'ark Po\'or]{M\'ark Po\'or$^\ast$}
\address{E\"otv\"os Lor\'and
University, Institute of Mathematics, P\'azm\'any P\'eter s. 1/c,
1117 Budapest, Hungary}
\email{sokmark@caesar.elte.hu}

\thanks{$^\ast$The first author was supported by the National Research, Development and Innovation Office
-- NKFIH, grants no. 124749, 129211}
\author[Saharon Shelah]{Saharon Shelah$^\dag$}
\address{Einstein Institute of Mathematics\\
	Edmond J. Safra Campus, Givat Ram\\
	The Hebrew University of Jerusalem\\
	Jerusalem, 91904, Israel\\
	and \\
	Department of Mathematics\\
	Hill Center - Busch Campus \\ 
	Rutgers, The State University of New Jersey \\
	110 Frelinghuysen Road \\
	Piscataway, NJ 08854-8019 USA}
\email{shelah@math.huji.ac.il}
\urladdr{http://shelah.logic.at}
\thanks{$^\dag$The second author was partially supported by the European Research
	Council (ERC) grant 338821, and by the Israel Science Foundation grant
	1838/19. Paper 1189 on Shelah's list.
 	\parbox[c]{11cm}{\scshape Supported by the \'UNKP-19-3 New National Excellence Program of the Ministry of Human Capacities.}
}


\title{Characterizing the spectra of cardinalities of branches of Kurepa trees}

\begin{abstract}
We give a complete characterization of the sets of cardinals that in a suitable forcing extension can be the  Kurepa spectrum, that is, the set of cardinalities of branches of Kurepa trees. 
This answers a question of the first named author.
\end{abstract}

\maketitle

\section{Introduction}

A tree is a Kurepa tree if it is of height $\omega_1$, each of its levels is countable, and it has more than $\omega_1$-many cofinal (that is of order type $\omega_1$) branches. In this paper we study the possible values of the branch spectrum of Kurepa trees, i.e. the set
\[ \Sp_{\omega_1} = \{ \lambda: \  \text{there exists a Kurepa tree } T \text{ s.t. } |\mathcal{B}(T)|=\lambda \} \ \subseteq [\omega_2, 2^{\omega_1}] \]
(where $\mathcal{B}(T)$ stands for the set of cofinal branches of $T$).

The spectrum is related to the model theoretical spectrum of maximal models of $\mathcal{L}_{\omega_1,\omega}$-sentences  \cite{SISO}.
Also canonical topological and combinatorial structures are associated with branches of Kurepa trees possessing a remarkably wide range
of nonreflecting properties \cite{Kos}. For higher  Kurepa trees (of weakly compact height) the consistency strength of certain types of the branch spectrum was studied in \cite{HM19}. 

It was first shown by Silver that the Kurepa Hypothesis (i.e. the existence of a Kurepa tree) is independent 
\cite{Si67}, or see \cite[Ch VIII, 3.]{Ku83}. 
Moreover the non-existence of Kurepa trees is equiconsistent with the existence of an inaccessible cardinal \cite[Ch VII, Ex. B8.]{Ku83}.

Questions about the possible values of the spectrum were addressed by Jin and Shelah in \cite{JiSh:469}. They proved  
(assuming an inaccessible cardinal) that consistently there are only Kurepa trees with $\omega_3$-many cofinal branches while $2^{\omega_1} = \omega_4$.

Building on ideas of Jin and Shelah,  the first named author provided a sufficient condition for a set to be equal to $\Sp_{\omega_1}$ in a forcing extension in \cite{PM}. Formally, it was shown that if $\mathbf{GCH}$ holds, and $0,1 \notin S$ is a set of ordinals such that $S$ satisfies either
	
	$	\underline{\text{Case A:}}$
\begin{enumerate}[label = (\roman*)]
	\item  $2 \in S$,
	\item \label{kettoo}  $\{\sup C: \ C \in [S]^{\leq \omega_1} \} \subseteq S$,
	\item $ (\forall \alpha \in S): \ \   (\omega\leq \cf(\alpha) < \omega_2) \ \to \  (\alpha+1 \in S) $,
\end{enumerate}	
	or
		 	
$\underline{\text{Case B:}}$ 
\begin{enumerate}[label = (\roman*)]
	\item $\exists \text{ an inaccessible } \kappa$,
	\item \label{Bkettoo} $\{ \sup C: \ C \in [S]^{<\kappa} \} \subseteq S$,
	\item $(\forall \alpha \in S): \ \  (\omega \leq \cf(\alpha) < \kappa) \ \to \  (\alpha+1 \in S)$,
\end{enumerate}	
then in a  forcing extension  we have $\{ \alpha: \ \aleph_\alpha \in \Sp_{\omega_1} \} = S$ (cardinals are only collapsed in Case B, from $(\omega_1, \kappa)$).
It can be easily seen that if $\cf(\mu)= \omega$ and $(\Sp_{\omega_1} \cap \mu)$ is cofinal in $\mu$, then there exists a Kurepa tree with $\mu$-many branches, as the union of countably many Kurepa trees is a Kurepa tree, and it is not difficult to see that the same holds if $\cf(\mu) = \omega_1$, therefore Case A / $\ref{kettoo}$ and Case B / $\eqref{Bkettoo}$ are in fact necessary. However, it remained a question whether the last clauses can be dropped.

In this paper as the main result we prove that assuming $\mathbf{CH} + (2^{\omega_1} = \omega_2)$ conditions (i), (ii) (in both cases)  are in fact sufficient by forcing a model of $\{ \alpha: \ \aleph_\alpha \in \Sp_{\omega_1} \} = S$. Also, we can arbitrarily prescribe $2^{\omega_1}$ to be any cardinal $\lambda \geq \sup (\Sp_{\omega_1})$ if in Case A  the equality $\lambda^{<\omega_2} = \lambda$ holds, or in Case B $\lambda^{<\kappa} = \lambda$ holds too. 

Moreover, when we do not want Kurepa trees with $\omega_2$-many cofinal branches, we prove that the inaccessible is necessary by verifying that if $\omega_2$ is a successor in $L$, then there exists a Kurepa tree with only $\omega_2$-many cofinal branches in $V$. It was known that these assumptions imply that there exists a Kurepa tree even in $L[A]$ for some $A \subseteq \omega_1$ \cite[Ch VII, Ex. B8.]{Ku83} (possibly having more than $\omega_2$-many cofinal branches in $V$). Our proof not only utilizes countable elementary submodels of initial segments of $L[A]$, but the nodes of the tree are such elementary submodels, and each cofinal branch uniquely corresponds to an initial segment of $L[A]$.

\section{Preliminaries, notations}
Under ordinals we always mean Neumann ordinals.
For a fixed cardinal $\chi$ we  will use the notation $\mathcal{H}(\chi)$ for the collection of sets of hereditary size less than $\chi$, i.e.
\[ \mathcal{H}(\chi) = \{ x: \ |\trcl(x)|  < \chi \}, \]
where $\trcl(x)$ stands for the transitive closure of $x$. In terms of forcing we will use the notations of \cite{Ku2013}, e.g. $p \leq q$ means that $p$ is the stronger. If it is clear from the context and won't make any confusion we will identify the set $x$ in the ground model with its canonical name $\check{x}$.
For a set $A$ the symbol $\mathcal{P}(A)$ denotes the powerset of $A$, and $[A]^{\lambda}$ stands for $\{X \in \mathcal{P}(A): \ |X| = \lambda \}$. 
For a function $s = \{ \langle \beta, s(\beta)\rangle: \ \beta \in \dom(s) \}$ we will also use the following notation and refer to $s$ as
\[ \langle s_\beta: \ \beta \in \dom(s) \rangle. \]  

Under a sequence we mean a function defined on a set of ordinals.
For sequences $s,t$  the relation $s = t \up \dom(s)$ (or equivalently $s \subseteq t$) will be also denoted by $s \vartriangleleft t$.

\bd A tree $\langle T, \prec_T \rangle$ is a partially ordered set (poset) in which for each $x \in T$
	the set
	\[ T_{\prec x} = \{ y \in T: \ y \prec_T x \} \]
	is well ordered by $\prec_T$.
\ed

\bd
	The height of $x$ in the tree $T$  is the  order type of $T_{\prec x}$
	\[ \htt(x) =  \otp(T_{\prec x}). \]
\ed

\bd
	For each ordinal $\alpha$ 
	the restriction of $T$ to $\alpha$  is 
	\[ T_{<\alpha} = \{ t \in T : \ \htt(t) < \alpha \} .\]
\ed

\bd
	The height of the tree $T$ (in symbols $\htt(T)$),  is the least $\beta$ such that
	\[ \nexists t \in T: \ \htt(t) = \beta.  \]
\ed

We will need the following lemma \cite[Ch II. Thm. 1.6.]{Ku83} which we will refer to as the $\Delta$-system Lemma.
\bl \label{delta}
	Let $\kappa$ be an infinite cardinal, let $\theta > \kappa$ be regular, and satisfy $\forall \alpha < \theta$ ($|\alpha^{<\kappa}| < \theta$). Assume that $|\mathcal{A}| \geq \theta$, and 
	$\forall x \in \mathcal{A}$ ($|x| < \kappa$). Then there is a $\mathcal{D} \subseteq \mathcal{A}$,
	such that $|\mathcal{D}|= \theta$, and $\mathcal{D}$ forms a $\Delta$-system, i.e.
	there is a kernel set $y$ such that 
	\[ \forall x \neq x' \in \mathcal{D}: \ \ x \cap x' = y. \]
\el

\section{The forcing} \label{forsz}
Now we can state our main theorem.
\bt \label{fot} Let $S_\bullet$ be a set of infinite cardinals such that $\omega, \omega_1 \notin S_\bullet$. Assume $\mathbf{CH}$, and that either 

$\underline{\text{Case 1:}}$
\begin{enumerate}[label = (\roman*)]
	\item $\omega_2 \in S_\bullet$,
	\item $2^{\omega_1} = \omega_2$,
	\item \label{harom} $\{ \sup C : \ C \in [S_\bullet]^{<\omega_2} \} \subseteq S_\bullet$,
\end{enumerate}
or

$\underline{\text{Case 2:}}$
\begin{enumerate}[label = (\roman*)]
	\item there exists an inaccessible $\kappa$ such that $S_\bullet \cap (\omega_1, \kappa) = \emptyset$,
	\item \label{harom'} $\{ \sup C : \ C \in [S_\bullet]^{<\kappa} \} \subseteq S_\bullet$.
\end{enumerate}

Then there exists a forcing extension $V^{\mathbb{P}}$ such that
\[ V^{\mathbb{P}} \models \ S_\bullet  = \Sp_{\omega_1}, \ \text{ where } \mathbb{P}  \text{ only collapses cardinals in } (\omega_1, \kappa) \text{ in } \underline{\text{Case2}}. \]
\et
The key will be Lemma $\ref{fol}$. After Lemma $\ref{lemma}$ we will put together the pieces in a short argument. Before these we need some preparation. 
\bd \label{kadf} In Case 1 (i.e. $\omega_2 \in S_\bullet$) define the cardinal $\kappa$ to be $\omega_2$.
\ed
\bcor No cardinal $\mu \notin (\omega_1, \kappa)$ is collapsed.
\ecor
\bt \label{fot2} Suppose that all conditions from Theorem $\ref{fot}$ hold, and $\kappa$ is defined in Definition $\ref{kadf}$. Assume further that $\lambda$ is a cardinal which is an upper bound of $S_\bullet$ such that $\lambda^{<\kappa} = \lambda$ (thus $\cf(\lambda) \geq \kappa$). Then there exists a forcing extension $V^\mathbb{P}$ with
\[ V^{\mathbb{P}} \models \ (S_\bullet = \{\mu: \ \text{there exists a Kurepa tree } T \text{ s.t. } |\mathcal{B}(T)| = \mu \}) \ \wedge \ (2^{\omega_1} = \lambda). \]
\et

\bd Let $S_\bullet^+ = S_\bullet \cup \{\kappa, \lambda \}$.
\ed


\bd \label{Pdf}
 For a cardinal $\theta \in S_\bullet$ let $\mathbb{Q}_\theta$ be the following notion of forcing.
 The triplet $p = \langle T_p, u_p, \overline{\eta}_p \rangle $ is an element of $\mathbb{Q}_\theta$ iff
 \begin{enumerate}[label=(\alph*)]
  	\item \label{PdfT} $T_p$ is a countable tree of height $\delta$ for some $\delta < \omega_1$ on the underlying set $\omega \cdot \delta$, where the $\beta$'th level is $[\omega \cdot \beta, \omega \cdot (\beta + 1))$, i.e. $T_{p, \leq \beta} \sm T_{p, < \beta} = [\omega \cdot \beta, \omega \cdot (\beta +1)) $ for each $\beta < \delta$,
 	\item \label{norm} for each $t \in T_p$ and $\beta < \delta$ there exists $t' \in T_p \sm T_{p,<\beta}$ s.t. $t \prec_{T_p} t'$,
 	\item $u_p \in [\theta]^{\leq \omega}$,
 	\item $\overline{\eta}_p = \langle \eta_{p,\alpha}: \ \alpha \in u_p \rangle$, where $\eta_{p,\alpha} \subseteq T_p$ is a branch in $T_{p,<\gamma}$ for some $\gamma \in \{\beta+1: \ \beta < \delta = \htt(T_p) \}$ (we do it for a technical reason,  we also could have stored only the maximal element instead of a chain with a maximal element).
 \end{enumerate}
Then $\mathbb{Q}_\theta$ is a poset with the obvious order, i.e. $q \leq p$, if $T_q$ is an end-extension of $T_p$, formally $T_{q, <\htt(T_p)} = T_p$, and for each $\alpha \in u_p$ the inclusion $\eta_{p,\alpha} \subseteq \eta_{q,\alpha}$ holds.
 
	Let $\name{T}_\theta, \name{\overline{\eta}}_\theta$ be the names for the generic tree and sequence, i.e. denoting  the generic filter by $\mathbf{G}_\theta$
	\[ \begin{array}{rll} 1_{\mathbb{Q}_\theta} \Vdash  & \name{T}_\theta = \cup\{ T_p: \ p \in \mathbf{G}_\theta \} &\text{ and} \\
 1_{\mathbb{Q}_\theta} \Vdash  & \name{\overline{\eta}}_\theta = \left\langle \name{\eta}_{\theta,\alpha} = \cup\{ \eta_{p,\alpha}: \ p \in \mathbf{G}_\theta \}: \ \ \alpha \in \theta \right\rangle. & \end{array} \]
\ed

\bd \label{Q*}
For a cardinal $\theta \in S_\bullet$ let $\mathbb{Q}^*_\theta \subseteq \mathbb{Q}_\theta$ be the following subposet.
\[ p \in \mathbb{Q}_\theta^*, \ \text { iff } \htt(T_p) \text{ is a successor, and } (\forall \alpha \in u_p): \ \eta_{p,\alpha} \text{ is a branch through } T_{p}.\]
\ed

\bd
If $\lambda \notin S_\bullet$ then let $\mathbb{Q}_\lambda$ be the countable supported product of $\langle ^{<\omega_1}2, \vartriangleleft \rangle$-s of length $\lambda$, i.e.
\[ \mathbb{Q}_\lambda = \{ p = \langle \eta_\alpha: \alpha \in u_p \rangle: \ (\forall \alpha \in u_p) \ \eta_\alpha \in \ ^{<\omega_1}2, \ \text{ for some } u_p \in [\lambda]^{\leq \omega} \}.  \]
\ed

\bd \label{Qkappa} If $\kappa \notin S_\bullet$ (and then $\kappa>\omega_2^V$ is inaccessible), then let $\mathbb{Q}_\kappa$ be the countable supported product of $\langle ^{<\omega_1}\gamma, \vartriangleleft \rangle$'s ($\gamma<\kappa$), a forcing which collapses each cardinal in $(\omega_1, \kappa)$:
\[ \mathbb{Q}_\kappa = \{ p = \langle \eta_\alpha: \ \alpha \in u_p \rangle: \ (\forall \alpha \in u_p) \ \eta_\alpha \in \ ^{<\omega_1}\alpha, \text{ for some } u_p \in [\kappa]^{\leq \omega}  \}. \]
\ed

\bd \label{Pkdf} We define the posets which we will need later.
\begin{enumerate}[label=\arabic*)]
	\item \label{PS} For $S \subseteq S^+_\bullet$ let $\mathbb{P}_S$ be the countable supported product of $\mathbb{Q}_\theta$-s ($\theta \in S$), i.e.
	\[ \mathbb{P}_S = \{ p \text{ is a function}: \ \dom(p) \in [S]^{\leq \omega} \ \wedge \ (\forall \theta \in \dom(p) \ p(\theta) \in \mathbb{Q}_\theta) \}.
	\]
	With a slight abuse of notation for $p \in \mathbb{P}_S$ and $\theta \in S \setminus \dom(p)$ we will mean $1_{\mathbb{Q}_\theta}$ under $p(\theta)$.
	\item For $\theta \in S_\bullet^+$, $U \subseteq \theta$ define its restriction from $\theta$ to $U$, i.e.
	\[ \mathbb{Q}_{\theta,U} = \{p \in \mathbb{Q}: \ u_p \subseteq U \}. \]
	
	\item For $S \subseteq S_\bullet^+$, $\overline{U} = \langle U_\theta: \theta \in S \rangle \in \prod_{\theta \in S} \mathcal{P}(\theta)$
	we define $\mathbb{P}_{S,\overline{U}}$ to be $\mathbb{P}$-s restriction to coordinates in $U_\theta$-s, i.e.
	\[ \mathbb{P}_{S,\overline{U}} = \{ p \in \mathbb{P}_S: \ (\forall \theta \in S) \ p(\theta) \in \mathbb{Q}_{\theta,U_\theta} \}. \]
	\item For $S, S' \subseteq S_\bullet^+$, $\overline{U} = \langle U_\theta: \theta \in S \rangle \in \prod_{\theta \in S} \mathcal{P}(\theta)$,  $\overline{U}' = \langle U'_\theta: \theta \in S \rangle \in \prod_{\theta \in S'} \mathcal{P}(\theta)$ we define
	\begin{itemize}
		\item $\overline{U} + \overline{U}' = \langle U_\theta \cup U'_\theta: \ \theta \in S \cup S' \rangle$ (where for $\theta \in S' \sm S$ under $U_\theta$ we mean the empty set, similarly for  $\theta \in S  \sm S'$, $U_\theta'$),
		\item $\overline{U} - \overline{U}' = \langle U_\theta \sm U'_\theta: \ \theta \in S \rangle$ (here we also mean the empty set  under $U'_\theta$ if $\theta \in S \sm S'$),
		\item $\overline{\id}_S = \langle \theta: \ \theta \in S \rangle$
		\item for the set $X$ if $\overline{W}_\alpha \in \prod_{\theta \in S} \mathcal{P}(\theta)$ ($\alpha \in X$) then  
		\[ \sum_{\alpha \in X} \overline{W}_\alpha = \left\langle \bigcup_{\alpha \in X} (W_{\alpha})_\theta: \ \theta \in S \right\rangle. \]
	\end{itemize}
\item Let $\mathbb{P} = \mathbb{P}_{S_\bullet^+}$.
	\item  If $p_0,p_0, \dots, p_n \in \mathbb{P}$ let $\bigwedge_{i\leq n } p_i$ denote the greatest lower bound if exists.
	
	
	\item \label{res} For $p \in \mathbb{P}$, and $S \subseteq  S_\bullet^+$, $\overline{U} = \langle U_\theta: \theta \in S \rangle \in \prod_{\theta \in S} \mathcal{P}(\theta)$ define $p \up \overline{U} \in \mathbb{P}_S$ to be the following restriction of $p \up S$ in the obvious fashion
	\[ \text{for each } \theta \in S: \ \ (p \up \overline{U})(\theta) = \langle T_{p(\theta)}, u_{p_\theta} \cap U_\theta, \overline{\eta}_p \up U_\theta \rangle. \]
\end{enumerate}
\ed

\bd \label{P*} For $S \subseteq S_\bullet^+$ define the notion of forcing $\mathbb{P}^*$ ($\mathbb{P}^*_{S}$, $\mathbb{P}^*_{S,\overline{U}}$, resp.) to be the subposet of $\mathbb{P}$ ($\mathbb{P}_{S}$, $\mathbb{P}_{S,\overline{U}}$, resp.) consisting of elements $p$ for that $p(\theta) \in \mathbb{Q}^*_\theta$ holds for each $\theta \in S_\bullet \cap \supp(p)$.

\ed

\br \label{P*P} The notion of forcing $\mathbb{P}^*$ ($\mathbb{P}^*_{S}$, $\mathbb{P}^*_{S,\overline{U}}$, resp.) is a dense subposet of $\mathbb{P}$ ($\mathbb{P}_{S}$, $\mathbb{P}_{S,\overline{U}}$, resp.), therefore forcing with $\mathbb{P}^*$ ($\mathbb{P}^*_{S}$, $\mathbb{P}^*_{S,\overline{U}}$, resp.) yields the same extensions as forcing with $\mathbb{P}$ ($\mathbb{P}_{S}$, $\mathbb{P}_{S,\overline{U}}$, resp.).
\er

\bc Let  $S \subseteq S^+_\bullet$, $\overline{U} = \langle U_\theta: \theta \in S \rangle$ be fixed. Then the poset $\mathbb{P}_{S,\overline{U}}$ has the $\kappa$-cc property.
\ec
\bp
Suppose that $\{p_\alpha: \ \alpha \in \kappa \} \subseteq \mathbb{P}_{S,\overline{U}}$ is an antichain. 
Working in $V'$, applying the $\Delta$-system lemma (Lemma $\ref{delta}$) for the system $\{ \dom(p_\alpha): \ \alpha \in \kappa \}$ of countable sets ($\ref{PS}$ from Definition $\ref{Pkdf}$), we obtain a set $A \in [\kappa]^\kappa$, such that the $\dom(p_\alpha)$'s ($\alpha \in A$) form a $\Delta$-system with kernel $K \subseteq S$. Since $K$ is obviously countable, for each $\alpha$ we have that $\langle T_{p_\alpha(\theta)}: \ \theta \in K \rangle$ is a countable sequence of countable trees (by $\ref{PdfT}$ from Definition $\ref{Pdf}$). This means that by $\mathbf{CH}$ we can assume that
\begin{equation} \label{uafa}  \langle T_{p_\alpha(\theta)}: \ \theta \in K \rangle = \langle T_{p_\beta(\theta)}: \ \theta \in K \rangle \ \ (\forall \alpha, \beta \in A). \end{equation}
 Now applying the $\Delta$-system lemma again for the system 
\[  U_\alpha = \bigcup_{\theta \in S} \left(\{\theta\} \times u_{p_\alpha(\theta)}\right)  \ \ (\alpha \in \kappa)  \]
yields a set $A' \in [A]^\kappa$ such that the $U_\alpha$'s ($\alpha \in A'$) form a $\Delta$-system with kernel $I \subseteq \bigcup_{\theta \in S} \{\theta\} \times \theta$ (of course, in fact, $I \subseteq \bigcup_{\theta \in K} \{\theta\} \times \theta$).
Now by $\eqref{uafa}$ it suffices to prove that 
\begin{equation} 
 \exists \alpha \neq \beta \in A' \ \text{ such that } (\text{for each } \langle \theta,\delta \rangle \in I): \ \ \eta_{p_\alpha(\theta), \gamma} =  \eta_{p_\beta(\theta), \gamma},
 \end{equation}
 for which it is enough to prove
 \beeq \label{kell}
 \left| \left\{ \left\langle \eta_{p_\alpha(\theta), \gamma}: \ \langle \theta, \gamma \rangle \in I \right\rangle: \ \alpha \in A' \right\} \right| < \kappa.
 \eeq
Fix $\alpha \in A'$. Now for each  $\langle \theta, \gamma \rangle \in I$, if  $\theta \in S_\bullet$  then $\eta_{p_\alpha(\theta), \gamma} \in \ [\omega_1]^{<\omega_1}$ (a branch through $T_{p_\alpha(\theta)}$). 

This means that (using that $I$ is countable) 
\begin{equation} \label{besz1} \left\{ \left\langle \eta_{p_\alpha(\theta), \gamma}: \ \langle \theta, \gamma \rangle \in I, \ \theta \in S_\bullet \right\rangle: \ \alpha \in A' \right\} \subseteq \prod_{\langle \theta, \gamma \rangle \in I, \ \theta \in S_\bullet} \  [\omega_1]^{<\omega_1}, \end{equation} 
which latter set is of size $\omega_1$ by $\mathbf{CH}$.
Second, if $\theta = \lambda \in (S_\bullet^+ \setminus S_\bullet) \cap S$, then 
\[ \left\{ \left\langle \eta_{p_\alpha(\theta), \gamma}: \ \langle \theta, \gamma \rangle \in I, \ \theta = \lambda \right\rangle: \ \alpha \in A' \right\} \subseteq \prod_{\langle \theta, \gamma \rangle \in I, \ \theta = \lambda} \ ^{<\omega_1}2. \]
Finally we have to consider the coordinate $\theta  = \kappa$ if $\kappa \in S \sm S_\bullet$. Then 
 letting $\delta = \sup\{\gamma: \ \langle \kappa, \gamma \rangle \in I \}$ we have $\delta < \kappa$, because $I$ is countable and $\kappa$ is inaccessible.
Then
\begin{equation} \label{besz2} \{ \langle \eta_{p_\alpha(\kappa), \gamma}: \ \langle \kappa, \gamma \rangle \in I \} \subseteq \prod_{\langle \kappa, \gamma \rangle \in I} \  ^{<\omega_1}\delta,
\end{equation}
and since $\kappa$ is inaccessible, this case $|\prod_{\langle \kappa, \gamma \rangle \in I} \  ^{<\omega_1}\delta| < \kappa$.
We obtain (using $\omega_1< \kappa$) that 
\[  \left| \{ \langle \eta_{p_\alpha(\theta), \gamma}: \ \langle \theta, \gamma \rangle \in I \} \right| \leq \omega_1 \cdot \omega_1 \cdot \left| \prod_{\langle \kappa, \gamma \rangle \in I} \  ^{<\omega_1}\delta \right| < \kappa,  \]
therefore $\eqref{kell}$ holds.
\ep

Now we make the intuition behind the easy idea of first adding the trees and some branches, and then forcing over the extension precise.
\bc \label{compemb} For each $S \subseteq S_\bullet^+$, $\overline{U} = \langle U_\theta: \ \theta \in S \rangle$ we have  \[ \mathbb{P}_{S,\overline{U}} \lessdot \mathbb{P}_S \lessdot \mathbb{P}, \] i.e. $\mathbb{P}_{S,\overline{U}}$ completely embeds into $\mathbb{P}_{S}$, which completely embeds into $\mathbb{P}$.
\ec
\bp Since $\mathbb{P} \simeq \mathbb{P}_S \times \mathbb{P}_{S_\bullet^+ \sm S}$, it is enough to prove that $\mathbb{P}_{S,\overline{U}} \lessdot \mathbb{P}_S$. 

Assume that $A \subseteq \mathbb{P}_{S,\overline{U}}$ is a maximal antichain in $\mathbb{P}_{S,\overline{U}}$, and let $p \in \mathbb{P}_S \sm \mathbb{P}_{S,\overline{U}}$. 
Then there exists $a \in A$, $a' \in \mathbb{P}_{S,\overline{U}}$ such that $a' \leq a$, $a' \leq b\up \overline{U}$. But then it is straightforward to check that also $a'$ and $b$ have a common lower bound.
\ep

\bd \label{Qc} Let $S \subseteq S_\bullet$, $\overline{U} = \langle U_\theta: \ \theta \in S \rangle$, $\theta_0 \in S$, $U'_{\theta_0} \subseteq \theta_0 \sm U_{\theta_0}$. Then $\name{\mathbb{Q}}^\circ_{\theta_0, U'_{\theta_0}} = \name{\mathbb{Q}}^\circ_{(S,\overline{U}), \theta_0, U'_{\theta_0}} $ denotes the $\mathbb{P}_{S,\overline{U}}$-name for a notion of forcing which adds the  branches $\name{\eta}_{\theta_0,\alpha}$  ($\alpha \in U'_{\theta_0}$) to $\name{T_{\theta_0}}$ in the following way
 \[  1 \Vdash_{\mathbb{P}_{S,\overline{U}}}    \name{\mathbb{Q}}^\circ_{\theta_0,U'_{\theta_0}} = \left\{ 
 \begin{array}{l} p = \langle \overline{\eta}_p, u_p \rangle:  \ (u_p \in [U'_{\theta_0}]^{\leq \omega}) \wedge  (\overline{\eta}_p = \langle \eta_{p,\alpha}: \ \alpha \in u_p \rangle) \text{,}  \\
 \text{such that each } \eta_{p,\alpha}  \text{ is a branch of } \name{T}_{\theta_0,< \delta_\alpha} \\ \text{for some } \delta_\alpha \in \{ \gamma+1: \ \gamma < \omega_1 \}  \\
    \end{array} \right\}. 
  \]
  If it is clear from the context we will use $\name{\mathbb{Q}}^\circ_{\theta_0, U'_{\theta_0}}$ not mentioning $S$ and $\overline{U}$.
\ed

\bd Let $S \subseteq S_\bullet$, $\overline{U} = \langle U_\theta: \ \theta \in S \rangle$, $\theta_0 \in S$.

If $\theta \in S_\bullet^+ \sm S_\bullet$, and $U'_\theta \subseteq \theta \sm U_\theta$, then define the $\mathbb{P}_{S,\overline{U}}$-name $\name{\mathbb{Q}}_{\theta,U'_{\theta}} = \name{\mathbb{Q}}^\circ_{\theta,U'_{\theta}}$ to be the name for $\mathbb{Q}_{\theta, U'_\theta}$.


\ed

\bd \label{DPh}
Let $S \subseteq S_\bullet^+$, $\overline{U} = \langle U_\theta: \ \theta \in S \rangle$, $\overline{U}' = \langle U'_\theta: \ \theta \in S \rangle \in \prod_{\theta \in S} \mathcal{P}(\theta)$, where $U_{\theta} \cap U'_{\theta} = \emptyset$ for each $\theta \in S$.
Then $\name{\mathbb{P}}^\circ_{\overline{U}'}=\name{\mathbb{P}}^\circ_{(S, \overline{U}), \overline{U}'}$  denotes the $\mathbb{P}_{S,\overline{U}}$-name 
for the countably supported product of $\name{\mathbb{Q}}^\circ_{\theta, U'_\theta}$'s ($\theta \in S$),
i.e. a notion of forcing which adds the branches $\name{\eta}_{\theta,\alpha}$  ($\alpha \in U'_{\theta}$) to $\name{T}_\theta$ for each $\theta \in S \sm S_\bullet$, and the sequences
$\name{\eta}_{\kappa,\alpha}$ ($\alpha \in U'_\kappa$) if $\kappa \in S \sm S_\bullet$, $\name{\eta}_{\lambda,\alpha}$ ($\alpha \in U'_\lambda$) if $\lambda \in S \sm S_\bullet$: 
\[   1 \Vdash_{\mathbb{P}_{S,\overline{U}}}    \name{\mathbb{P}}^\circ_{\overline{U}'} = \left\{ p \text{ is a function}: \ \dom(p) \in [S]^{\leq \omega} \ \wedge \ (\forall \theta \in \dom(p) \ p(\theta) \in \name{\mathbb{Q}}^\circ_{\theta,U'_{\theta}}) \} \right\}. \]
Again, as in Definition $\ref{Qc}$ if it does not cause any confusion we only use the notation $\name{\mathbb{P}}^\circ_{\overline{U}'}$ not mentioning $S$ and $\overline{U}$.
\ed

The following claim is an easy observation. 


\bc \label{rc2}
If $\mathbf{G}$ is a  $\mathbb{P}_{S,\overline{U}}$-generic filter over $V$ (where $S \subseteq S_\bullet^+$, $\overline{U} = \langle U_\theta: \ \theta \in S \rangle$, $\overline{U}' = \langle U'_\theta: \ \theta \in S \rangle \in \prod_{\theta \in S} \mathcal{P}(\theta)$, and $U_{\theta} \cap U'_{\theta} = \emptyset$ for each $\theta \in S$), then 
with the notation from \cite{Ku2013} 
\[ \mathbb{P}_{S,\overline{U} + \overline{U}'} / \mathbf{G} = \{ p \in \mathbb{P}_{S,\overline{U} + \overline{U}'}: \ \forall q \in \mathbf{G} \ p \not\perp q \},\]
the quotient poset $\mathbb{P}_{S,\overline{U} + \overline{U}'} / \mathbf{G}$ and the evaluation of $\name{\mathbb{P}}^\circ_{\overline{U'}}$ are isomorphic, i.e.
	\[V[\mathbf{G}] \models \ \ \name{\mathbb{P}}^\circ_{\overline{U}'}[\mathbf{G}] \simeq \mathbb{P}_{S,\overline{U} + \overline{U}'} / \mathbf{G}. \]
\ec

Since $\mathbb{P}_{S,\overline{U}}$ completely embeds into $\mathbb{P}_{S,\overline{U} + \overline{U}'}$ (by Claim $\ref{compemb}$),
\cite{Ku2013}[Lemma V.4.45.] (and \cite[Lemma V.4.44.]{Ku2013}) implies the following.

\bc \label{canon}
 Let $S \subseteq S_\bullet^+$, $\overline{U} = \langle U_\theta: \ \theta \in S \rangle$, $\overline{U}' = \langle U'_\theta: \ \theta \in S \rangle \in \prod_{\theta \in S} \mathcal{P}(\theta)$, where $U_{\theta} \cap U'_{\theta} = \emptyset$ for each $\theta \in S$.
 Then the canonical embedding from $\mathbb{P}_{S,\overline{U} + \overline{U}'}$ to the iteration $\mathbb{P}_{S,\overline{U}} \ast (\mathbb{P}_{S,\overline{U} + \overline{U}'}/\mathbf{G})$ is a dense embedding.
\ec


Now putting together Claims $\ref{rc2}$ and $\ref{canon}$ we have the following, meaning that instead of forcing with $\mathbb{P}_{S,\overline{U} + \overline{U}'}$ we can force with $\mathbb{P}_{S,\overline{U}}$ and then with (the evaluation of) $\name{\mathbb{P}}^\circ_{\overline{U}'}$.

\bl \label{mindegyforsz} Let $S \subseteq S_\bullet^+$, $\overline{U} = \langle U_\theta: \ \theta \in S \rangle$, $\overline{U}' = \langle U'_\theta: \ \theta \in S \rangle \in \prod_{\theta \in S} \mathcal{P}(\theta)$, where $U_{\theta} \cap U'_{\theta} = \emptyset$ for each $\theta \in S$.
Then forcing with $\mathbb{P}_{S,\overline{U} + \overline{U}'}$ amounts to the same as forcing with
$\mathbb{P}_{S,\overline{U}}$ and then with $\mathbb{P}_{S,\overline{U} + \overline{U}'} / \mathbf{G} \simeq \name{\mathbb{P}}^\circ_{\overline{U}'}$.

\el

\bd
If $S \subseteq S_\bullet^+$, $\overline{U} = \langle U_\theta: \ \theta \in S \rangle$, $\overline{U}' = \langle U'_\theta: \ \theta \in S \rangle \in \prod_{\theta \in S} \mathcal{P}(\theta)$. Now if $\mathbf{G}$ is generic over $\mathbb{P} = \mathbb{P}_{S_\bullet^+}$ then we define 
\begin{itemize}
 \item $\mathbf{G}_S = \mathbf{G} \cap \mathbb{P}_S$,
  \item  $\mathbf{G}_{S, \overline{U}} = \mathbf{G} \cap \mathbb{P}_{S, \overline{U}}$,
  \item and $\mathbf{G}^\circ_{\overline{U}'} \subseteq \mathbb{P}^\circ_{\overline{U}'}[G_{S,\overline{U}}] \in V[G_{S,\overline{U}}]$ to be the filter given by the canonical mapping from Claims $\ref{rc2}$, $\ref{canon}$.
\end{itemize}
\ed

The following are basic observations. Roughly speaking, we isolate a dense $\omega_1$-closed subset of a two-step iteration similarly as in \cite{Kunen78}.
\bc $\mathbb{P}^*$ (and in general each $\mathbb{P}^*_{S,\overline{U}}$) is $\omega_1$-closed, i.e. for each decreasing sequence of type $\omega$ has a lower bound. In particular if $\mathbf{G}^* \subseteq \mathbb{P}^*$, (or in general $\mathbf{G}^*_{S,\overline{U}} \subseteq \mathbb{P}^*_{S,\overline{U}}$) is generic over $V$, then there is no new sequence of ordinals of type $\omega$.
\ec

The last claim and Remark $\ref{P*P}$ obviously implies the following.
\bc \label{ome}
 Forcing with $\mathbb{P}$ (or $\mathbb{P}_{S,\overline{U}}$) doesn't add new sequence of ordinals of type $\omega$, and for a given generic filter $\mathbf{G} \subseteq \mathbb{P}$ 
 \[ \mathcal{H}(\omega_1)^V = \mathcal{H}(\omega_1)^{V[\mathbf{G}]} = \mathcal{H}^{V[\mathbf{G}_{S,\overline{U}}]}. \]
\ec

\bl \label{her} Let $\mathbf{G} \subseteq \mathbb{P}_{S, \overline{U}}$ generic over $V$, $B \in V[\mathbf{G}]$ where $B \subseteq \mathcal{H}(\omega_1)$. Then (in $V$) there exist $S_* \subseteq S$, $|S_*| < \kappa$ and $\overline{W}_* = \langle W^*_\gamma: \ \gamma \in S_* \rangle \in \prod_{\gamma \in S_*} [U_\gamma]^{<\kappa}$, such that $B \in V[\mathbf{G}_{S_*,\overline{W}_*}]$.
\el
\bp
Choose $p \in \mathbf{G}$ forcing that $B \subseteq \mathcal{H}(\omega_1)$, and a nice $\mathbb{P}_{S, \overline{U}}$-name for $B$, obtaining for each $x \in \mathcal{H}(\omega_1)$ an antichain $A_x \subseteq \mathbb{P}_{S, \overline{U}}$ deciding about $x \in B$. Then by $\kappa$-cc  we have that each $|A_x| < \kappa$, the set $S_\ast = \bigcup_{x \in \mathcal{H}(\omega_1)} \bigcup_{a \in A_x} \dom(a)$ is of size less than $\kappa$ (as $\kappa$ is either inaccessible, or $\omega_2$).
Also for $\theta \in S_*$ the set $W^*_\theta =  \bigcup_{x \in \mathcal{H}(\omega_1)} \bigcup_{a \in A_x} u_{a(\theta)}$ is smaller that $\kappa$. Now it is easy to see that $\overline{W}_* = \langle W^*_\gamma: \ \gamma \in S_* \rangle$ is as claimed.  
\ep

Then the following immediately follows from the $\omega_1$-closedness, and $\kappa$-cc.
\bc \label{card}
Forcing with $\mathbb{P}$ doesn't collapse $\omega_1$, and cardinals at least $\kappa$. Moreover, if $\mathbf{G} \subseteq \mathbb{P}$ is generic, then
\[ V[\mathbf{G}] \models  \ \text{"}\kappa = \omega_2 \text{"}. \]
\ec

\bl \label{fol}
Let $T \in V[\mathbf{G}_{S,\overline{U}_*}]$ be a Kurepa tree, $S' \subseteq S$ ($S' \in V$). 
Then, if $b \in V[\mathbf{G}_{S,\overline{U}_* + \id_{S'}}]$ is a branch of $T$, then there exists a finite set $S'' \subseteq S'$, and $\overline{U}_\bullet = \langle U^\bullet_\theta: \ \theta \in S'' \rangle$ s.t. each $U^\bullet_\theta$ is finite, and $b$ is in the model obtained by adding these finitely many $\eta_{\theta,\alpha}$'s ($\theta \in S''$, $ \alpha \in U^\bullet_\theta$) to $V[\mathbf{G}_{S,\overline{U}_*}]$, i.e. 
\[ b \in V[\mathbf{G}_{S,\overline{U}_* + \overline{U}_\bullet}]. \]

\el
\bp

Let $\dot{T} \in V$ be a $\mathbb{P}_{S,\overline{U}_*}$-name for $T$.
Define
\beeq
\mathbb{P}' = \mathbb{P}_{S,\overline{U}_* + \overline{\id}_{S'}}.
\eeq
Suppose that $p_* \in \mathbb{P}'$ forces that $\dot{b} \in V$ is a $\mathbb{P}'$-name for a counterexample (i.e. forcing that for no such $\overline{U}_\bullet$ there exists a $\mathbb{P}_{\overline{U}_* + \overline{U}_\bullet}$-name $\dot{b}'$ - which is of course also a $\mathbb{P}'$-name - with $\dot{b}' = \dot{b}$). Let $\chi$ be large enough, and let $\langle N_0, \in \rangle \prec \langle \mathcal{H}(\chi), \in \rangle$ be countable s.t. $p_*,\dot{b}, \dot{T}, S, S', \overline{V}, \mathbb{P}_{S,\overline{U}_*} \in N_0$.

Let $\delta_\bullet = N_0 \cap \omega_1$. Define the countable set $N_1$ to be such that $N_0 \in N_1$, and $\langle N_1, \in \rangle \prec \langle \mathcal{H}(\chi ), \in \rangle$. 
Let $X$ be set of the indices of the new branches  added to $\langle \name{T}_\theta: \ \theta \in S' \rangle$ by $\mathbf{G}_{S,\overline{U}_* + (\id_{S'})}$ that are in  $V[\mathbf{G}_{S,\overline{U}_* +\overline{\id}_{S'}}] \setminus V[\mathbf{G}_{S,\overline{U}_*}]$, and belong to $N_0$, i.e.
\beeq \label{Xdf} X= N_0 \cap \{\langle \theta, \alpha\rangle : \ (\theta \in S') \ \wedge \ (\alpha \in \theta \sm U^*_\theta) \}.\eeq
We fix an enumeration of $X$ and define also the sequence of the first $n$ indices from this countable set, and as well for each $n$ the one-length sequence consisting only the $n$'th, that is
\[
\begin{array}{c} \text{let } \ \langle \langle \varrho_n, \xi_n \rangle: \ n \in \omega, n>0 \rangle \text{ enumerate } \ X \text{ (starting from }1), \end{array} \]
\beeq \label{Wdf}
\begin{array}{rl}
	\overline{W}_n = & \langle W_{n,\theta}: \ \theta \in S' \cap N_0 \rangle, \\
	& \text{where } W_{n,\theta}= \{ \alpha: \ \langle \theta, \alpha \rangle= \langle \varrho_j, \xi_j \rangle  \text{ for some } j\leq n \} \\
	\overline{w}_n = & \langle w_{n,\theta}: \ \ \theta \in S' \cap N_0 \rangle  \\
	 & \text{ where } w_{n,\theta} = \{ \xi_{n} \} \ \text{ if } \theta = \varrho_{n}, \ w_{n,\theta} = \emptyset \ \text{ otherwise}.
\end{array}
\eeq

Observe that if $p \in \mathbb{P} \cap N_0$, then each $\theta \in \dom(p)$ is an element of $N_0$  since $\dom(p)$ is countable (by Definition $\ref{Pkdf}$), and similarly $T_{p(\theta)}, u_{p(\theta)} \subseteq N_0$ (by Definitions $\ref{Pdf}-\ref{Qkappa}$).

Working in $V$ we will construct an $N_0$-generic condition in $\mathbb{P}'$, which will derive us to a contradiction.
It is enough to prove the following claim.
\bc \label{focl}
There exists a sequence $\langle \overline{p}_n: \ n \in \omega \rangle \in V$, $p'_0 \in \mathbb{P}_{S, \overline{U}*}$ and a sequence $\overline{q} = \langle q_n: \ n \in \omega \rangle$ such that the following holds.
 \begin{enumerate}[label=$\boxplus_{\text{\arabic*}}$, ref=\arabic*]
 	\item \label{pv} $\overline{p}_0 = \langle p_{0,l}: \ l \in \omega \rangle$ is such that \begin{enumerate}[label=(\alph*)]
 			\item $p_{0,0} = p_* \up \overline{U}_ *$,
 			\item $p_{0,l} \in N_0 \cap \mathbb{P}_{S, \overline{U}_*}$ for each $l \in \omega$,
	 		\item $\langle p_{0,l}: \ l \in \omega \rangle$ is $\leq_{\mathbb{P}}$-decreasing,
	 		\item $\overline{p}_0 \in N_1$,
	 		\item letting $\mathbf{G}_0 = \{ p \in \mathbb{P}_{S,\overline{U}*} \cap N_0: \ (\exists l) \ p \geq p_{0,l} \}$, the filter $\mathbf{G}_0$ is $\mathbb{P}_{S,\overline{U}_*}$-generic over $N_0$.
 	\end{enumerate}
 	\item \label{p0} $p'_0\in \mathbb{P}_{S,\overline{U}*}$ satisfies the following 
 	  \begin{enumerate}[label=(\alph*)]
	 	  	\item $p'_0$ is a lower bound of $p_{0,l}$ for each $l \in \omega$ (hence forces a value to $\name{T}_{\theta,<\delta_\bullet}$ for each $\theta \in S \cap N_0$),
	 	  	\item \label{p00} $p'_0$ forces a value  to 	
	 	  	   $\name{T}_{\theta,\leq\delta_\bullet}$ for each $\theta \in S \cap N_0$ such that for every $\delta_\bullet$-branch $B$ in $\name{T}_{\theta,<\delta_\bullet}$ the inclusion $B \in N_1$ implies that $B$ has an upper bound in  $\name{T}_{\theta,\leq\delta_\bullet}$,
	 	  	\item \label{p000} $p'_0$ forces a value to $\dot{T}_{\leq\delta_\bullet}$.
 	  	
 	  	\end{enumerate}
   	\item \label{pvonas} for every $n >0 $ the sequence $\overline{p}_n = \langle p_{n,l}: \ l \in \omega \rangle$ has the following properties.
   		\begin{enumerate}[label=(\alph*)]
   			\item $\forall l \in \omega$  $p_{n,l} \in N_0 \cap \mathbb{P}_{S,\overline{U}_* + \overline{w}_{n}}$,
   			\item \label{pnG0} $p_{n,l} \upharpoonright \overline{U}_* \in \mathbf{G}_0$
			\item $\langle p_{n,l}: \ l \in \omega \rangle$ is $\leq_{\mathbb{P}}$-decreasing,
			
   			\item \label{N1ben} $\overline{p}_n \in N_1$,
   			\item \label{pngen} letting 
   			\[\mathbf{G}_n = \{  p \in \mathbb{P}_{S,\overline{U}* + \overline{W}_n} \cap N_0: \ (\exists l_0,l_1,\dots, l_n) \ p \geq \bigwedge_{j=0}^{n} p_{j,l_j} \}, \] the filter $\mathbf{G}_n$ is $\mathbb{P}_{S,\overline{U}_*+ \overline{W}_n}$-generic over $N_0$.
   		\end{enumerate}
	 \item \label{qv} For the sequence $\overline{q} = \langle q_n: \ n \in \omega \rangle$
	 \begin{enumerate}[label=(\alph*)]
		 	\item $q_n \in N_0 \cap \mathbb{P}_{S, \overline{U}_* + \overline{\id}_{S'}}$ for each $n \in \omega$,
		  \item $q_0 = p_*$,

		  \item \label{qsuc} $\langle q_n: \ n \in \omega \rangle$ is $\leq_{\mathbb{P}}$-decreasing,
		  \item \label{qgen} $\forall n$: $q_n \upharpoonright (\overline{U}_* + \overline{W}_n) \in \mathbf{G}_n$,
	 \item \label{qv1} 
	
	 Let $\langle \dot{B}_n: \ n \in \omega \rangle$ enumerate the branches of $\dot{T}_{<\delta_\bullet}$ which has an upper bound in $\dot{T}_{\leq \delta_\bullet}$ (forced by $p'_0$).
	 Then  $q_{n+1} \wedge p'_0$ forces that $\dot{b} \neq B_n$, which will be guaranteed by the following requirement: 
	 
	 There exist $\delta < \delta_{\bullet}$, $t \neq t' \in \dot{T}_{\leq \delta} \sm \dot{T}_{<\delta}$, such that $p'_0$ forces $B_n$-s $\delta$'th level to be $t'$, and  $q_{n+1}$ forces $t \in \dot{b}$, i.e. 
	  \beeq \label{nemBn} \begin{array}{l} p'_0 \Vdash \ \dot{B}_n \cap (\dot{T}_{\leq \delta} \sm \dot{T}_{<\delta}) = \{t'\}  \\
	  	\text{and} \\
   	 q_{n+1} \Vdash  \ \dot{b} \cap (\dot{T}_{\leq \delta} \sm \dot{T}_{<\delta}) = \{t \}. \end{array} \eeq
	(Observe that the latter is a statement in $N_0$.)
 	\end{enumerate}
 \end{enumerate}
\ec 
Before proving Claim $\ref{focl}$ we argue why this claim implies Lemma $\ref{fol}$. First, the claim gives the following condition in $\mathbb{P}_{S, \overline{U}_* + \overline{\id}_{S'}}$. For each $n \in \omega$ let $\eta_{\varrho_n,\xi_n}$ be the branch in $\name{T}_{\varrho_{n}, < \delta_\bullet}$ represented by the sequence $\overline{p}_{n}$, i.e.
\beeq \eta_{\varrho_{n},\xi_{n}} = \cup \{ \eta_{p_{n,l}(\varrho_{n}),\xi_{n}}: \ l \in \omega \}, \eeq
and note that $\eta_{\varrho_n,\xi_n} \in N_1$ ($n \in \omega$) by $\bref{pvonas}/\ref{N1ben}$. Therefore by $\bref{p0} / \ref{p00}$ we can extend each $\eta_{\varrho_n,\xi_n}$ to a branch $\eta'_{\varrho_n,\xi_n}$ in $(T_{p'_0(\varrho_n)})_{< \delta_\bullet +1}$.
Define the function $p_\bullet$ to be the extension of $p'_0$ by the $\eta_{\varrho_n, \xi_n}$'s in the obvious way: (Note that by $\bref{p0}$ we have $S \cap N_0 \subseteq \dom(p'_0) \subseteq S$, and for each $\theta \in S \cap N_0$ the inclusion $U^*_\theta \cap N_0 \subseteq u_{p_0'(\theta)} \subseteq U^*_\theta$.) Define $p_\bullet$ to be function on $\dom(p'_0)$ such that if $\theta \notin N_0 \cap S'$, then $p_\bullet(\theta) = p'_0(\theta)$, and for $\theta \in N_0 \cap S'$ define $p_\bullet(\theta)$ to be the following proper extension of $p'_0(\theta)$.
Let $u_{p_\bullet(\theta)} = u_{p_0(\theta)} \cup (\theta \cap N_0)$, and if $\alpha \notin  u_{p'_0(\theta)}$ (when necessarily $\alpha \notin U^*_\theta$) and (by $\eqref{Wdf}$) choose $n >0 $ so that   
\beeq \langle \theta, \alpha \rangle = \langle \varrho_{n}, \xi_{n} \rangle, \text{ and let } \eta_{p_\bullet(\theta),\alpha} = \eta'_{\varrho_n,\xi_n}, \eeq 
otherwise 
\beeq \eta_{p_\bullet(\theta),\alpha} = \eta_{p_0(\theta),\alpha} \ \ \ \ \ (\text{if }\alpha \in U^*_\theta). \eeq 

Observe that as $\eta'_{\varrho_n,\xi_n}$  was a cofinal branch in $(T_{p_\bullet(\varrho_n)})_{< \delta_\bullet+1}= (T_{p'_0(\varrho_n)})_{< \delta_\bullet+1}$ our function $p_\bullet$ is indeed a condition in $\mathbb{P}_{S,\overline{U}_* + \overline{\id}_{S'}}$. Moreover, the following shows that $\forall n \in \omega$ $p_\bullet \leq q_n$. Fix $n\in \omega$, then using $\bref{qv}/ \ref{qgen}$ we have $q_n \upharpoonright (\overline{U}_* + \overline{W}_n) \in \mathbf{G}_n$, i.e. there exist $l_0,l_1,\dots l_n \in \omega$, such that $\bigwedge_{j=0}^{n} p_{j,l_j} \leq_{\mathbb{P}} q_n \upharpoonright (\overline{U}_* + \overline{W}_n)$.
This means that 
\[\bigwedge_{j=0}^{n} p_{j,l_j} \leq q_n \upharpoonright (\overline{U}_* + \overline{W}_0) = q_n \upharpoonright (\overline{U}_*), \] and 
\[ \text{for each }0 < j\leq n \ \ \eta_{q_n(\varrho_{j}),\xi_{j}} \subseteq \eta_{p_{j,l_j}(\varrho_{j}), \xi_{j}} \subseteq \eta'_{\varrho_{j}, \xi_{j}} = \eta_{p_\bullet(\varrho_{j}), \xi_{j}} . \]
  On the other hand, for $j>n$ we have (recalling $\overline{q} = \langle q_n: \ n \in \omega \rangle$ is $\leq_{\mathbb{P}}$-decreasing by $\bref{qv}$)  that
\[ \eta_{q_n(\varrho_{j}),\xi_{j}} \subseteq \eta_{q_j(\varrho_{j}),\xi_{j}} \subseteq \eta'_{\varrho_{j},\xi_{j}} = \eta_{p_\bullet(\varrho_{j}), \xi_{j}},\] 
therefore $p_\bullet \leq q_n$, indeed.

Now assuming $p_\bullet \in \mathbf{G}_{S,\overline{U}_* + \overline{id}_{S'}}$ will easily yield a contradiction:
First recall that $p_*$ (and therefore as well $q_0$ and $p_\bullet$) forced that $\dot{b}$ is a branch through $\dot{T}$. Then $\bref{p0}/ \ref{p000}$ implies that $p'_0$, thus $p_\bullet$ as well determines $\dot{T}_{\leq \delta_\bullet}$, and $p_\bullet$ forces (by $\bref{qv} / \ref{qv1}$) that each element of the $\delta_\bullet$'th level of $\dot{T}$ is the upper bound of $B_i$ for some $i \in \omega$. This means that
\[ p_\bullet \Vdash \ (\exists i \in \omega) \  \dot{b} \cap \dot{T}_{<\delta_\bullet} = B_i,\]
while at the same time
\[ (q_i\wedge p'_0) \Vdash \dot{b} \neq B_i, \] 
since $\eqref{nemBn}$ holds.

This together with $p_\bullet \leq q_i,p'_0$ gives the contradiction.
Now we can turn to the proof of the claim.

\bp(Claim $\ref{focl}$)

For the construction of each sequence $\overline{p}_n$ and each $q_n$ we will work in $N_1$. 
This will need a lot of preparation.

Recall that $X \subseteq N_0$ denoted the indices of branches added by forcing with $\mathbb{P}_{S,\overline{U}_* + \overline{\id}_{S'}} \cap N_0$ but missing from $V[\mathbf{G}_{S,\overline{U}_*}]$ $\eqref{Xdf}$, and that for each condition $p$, $\theta \in S_\bullet$, and $\delta < \omega_1$ the $\delta$'th level of $T_{p(\theta)}$ is (a subset of) $[\omega \cdot \delta, \omega \cdot (\delta+1))$.
Define $E \subseteq N_0$ as follows.
\beeq \label{eDf} \begin{array}{rl} e \in E  \ \text{ iff } &     e \in N_0, \text { and }e =(u_e, \overline{\eta}_e), \text{ where } u_e \in [X]^{\leq \omega}, \\
& \overline{\eta}_e=\langle \eta_{e,\theta,\alpha}: \ \langle\theta,\alpha \rangle \in u_e \rangle, \ \text{ such that } \\
&  \eta_{e,\theta,\alpha} \subseteq \omega \cdot (\delta_{\theta,\alpha}+1) \text{ for some } \delta_{\theta, \alpha} < \omega_1  \\

\end{array} \eeq
\bd
For each $n$, $p \in \mathbb{P}_{S,\overline{U}_* + \overline{W}_n}$, and $e \in E$, \underline{\text{if}} for each $\langle \theta,\alpha \rangle \in u_e$ we have $\theta \in \dom(p)$, and for each $i<n$  $\langle \varrho_i, \xi_i \rangle \notin u_e$ holds then define $p \tieconcat e$ as
\beeq \label{pce} \begin{array}{l}
	\dom(p \tieconcat e) = \dom(p), \\
	u_{(p \tieconcat e)(\theta)} = u_{p(\theta)} \cup \{ \alpha: \ \langle \theta,\alpha \rangle \in u_e\} ) \ \ (\forall \theta \in \dom(p \tieconcat e) ), \\
	\eta_{(p \tieconcat e)(\theta), \alpha} = \left\{ \begin{array}{l} \eta_{p(\theta),\alpha}, \ \ \text{ if }  \alpha \in u_{p(\theta)}, \\
		\eta_{e,\theta,\alpha}, \ \ \ \text{ if } \langle \theta, \alpha \rangle \in u_e,  \end{array} \right. \\
	\text{if this is a condition in }\mathbb{P} \text{ (i.e. for each }\langle \theta,\alpha \rangle  \in u_{e} \\
	\eta_{e,\theta,\alpha} \text{ is a cofinal branch of } (T_{p(\theta)})_{<\delta+1} \text{ for some } \delta \leq \htt(T_{p(\theta)}), \\
	\\
	\text{ otherwise }p \tieconcat e = \emptyset.
	
\end{array} \eeq
\ed

Let $\mathcal{D}$ denote the set of dense subsets of $\mathbb{P}_{S,\overline{U}_* + \overline{\id}_{S'}}$. Fix an enumeration
\[ \left\langle \langle J_i, \varepsilon_i \rangle : \ i \in \omega \right\rangle \in N_1  \ \ \text{ of } \ (\mathcal{D} \cap N_0) \times E, \]
and let $k(D,e)$ denote the index of the pair $\langle D, e \rangle$ (i.e. 
\beeq \label{kdf} J_{k(D,e)} = D, \ \varepsilon_{k(D,e)} = e\text{), then we also have } k \in N_1, \text{ of course.}\eeq
 Fix a function $g \in N_0$
\beeq \label{gfv}
 \begin{array}{rl} 
	g: & \mathbb{P}_{S,\overline{U}_*+\overline{\id}_{S'}} \times \mathcal{D} \ \to \ \mathbb{P}_{S,\overline{U}_*+\overline{\id}_{S'}}  \\
	 \text{with } & \forall p,D:  \\ 
    &	\bullet_1 g(p,D) \in D, \\
	& \bullet_2 \ g(p,D) \leq p, \\
	
	 \end{array}
	 \eeq
\[ \text{(Then } g \in N_0 \text{ obviously implies }(p,D \in N_0 \ \Rightarrow \ g(p,D) \in N_0 \text{).)} \]
We will have to define also the auxiliary sequence $\overline{r} = \langle r_l : \ l \in \omega \rangle$ with the following property:
\begin{enumerate}[label=$\circledast_{\arabic*}$, ref=\arabic*]
	\item \label{p'1}$\overline{r} \in N_1$,
	\item for each $l$ $r_{l} \in \mathbb{P}_{S, \overline{U}_*} \cap N_0$,
	\item for each $l$ $p_{0,l+1} \leq r_{l} \leq p_{0,l}$,
	\item \label{p'l}if there exists $p \in \mathbb{P}_{S,\overline{U}^*}$ such that $p \leq p_{0,l}$, and $p\tieconcat \varepsilon_{l}$ is a condition  extending $p_{0,l}$ in  $\mathbb{P}_{S, \overline{U}^* + \overline{\id}_{S'}}$, then $r_{l}$ is such that.
\end{enumerate}

Now we can construct the $p_{0,i}$'s (and $r_i$'s). Let $p_{0,0} = p_* \up \overline{U}_*$.
For obtaining the $p_{0,l}$'s proceed as follows. Assume we have defined $p_{0,0},p_{0,1}, \dots ,p_{0,l-1}$ (and as well  the $r_i$'s for $i< l-1$). Now if there exists $p \in \mathbb{P}_{S,\overline{U}^*}$ $p \leq p_{0,l-1}$, s.t. $p \tieconcat \varepsilon_{l-1} \neq \emptyset$ but a condition extending $p_{0,l-1}$, then let $r_{l-1} \in N_0$ be such a $p$ (recall that $\varepsilon_{l-1} \in E \subseteq N_0$ by $\eqref{eDf}$), otherwise define $r_{l-1}= p_{0,l} = p_{0,l-1}$. Lastly, in the former case define $p_{0,l} = g(r_{l-1}, D_{l-1}) \up \overline{U}_*$. It is clear from the construction and the definition of $g$ that $p_{0,l-1} \leq r_{l-1} \leq p_{0,l}$, and $r_{l-1}, p_{0,l}\in N_0$, and since every object as well as the series $\langle \varepsilon_i: \ i \in \omega \rangle$ are elements of $N_1$, we obtain  $\overline{p}_0, \overline{r}_0 \in N_1$, too.

Finally, it is straightforward to check that the filter $\mathbf{G}_0$ generated by the $p_{0,l}$'s meets every dense subset $D \in N_0$ of $\mathbb{P}_{S,\overline{U}_*}$. Fixing such a $D$ 
\[ D' = \{p \in \mathbb{P}_{S,\overline{U}_*+ \overline{\id}_{S'}}: \ p \up \overline{U}_* \in D \} \] is clearly a dense subset of $\mathbb{P}_{S,\overline{U}_* + \overline{\id}_{S'}}$ belonging to $N_0$. 
This means that if $e \in E$ is the empty sequence, then there exists $i \in \omega$, such that $J_i= D'$, and $\varepsilon_i = e$, therefore $p_{0,i+1} \in D$.

For $p'_0$, first consider the condition $p''_0 \in N_1$ consisting of only the generic trees given by $\mathbf{G}_0$ (for each $\theta \in \dom(p''_0)=  N_0 \cap S$ the tree $T_{p'_1(\theta)} = \cup\{ T_{p(\theta)}: p \in \mathbf{G}_0 \}$ is of height $\delta_\bullet$, but $u_{p''_0(\theta) = \emptyset}$). Then let $p'''_0 \in \mathbb{P}_{S, \overline{U}_*}$, $p'''_0 \leq p''_0$ be an extension so that for each $\theta \in S' \cap N_0$ the tree $T_{p'_2}(\theta)$ satisfies that for each branch $B$ through $(T_{p'''_0}(\theta))_{<\delta_\bullet} = T_{p''_0(\theta)}$, if $B \in N_1$, then there is an upper bound of $B$ in $T_{p'''_0(\theta)}$. This can be done since $N_1$ is countable. Moreover, we choose the other part of  $p'''_0$ so that for each $\theta, \alpha \in N_0$, if $\alpha \in U^*_\theta$ the chain $\eta_{p'''_0(\theta),\alpha}$ (with a top element) contains the chain $\cup\{ \eta_{p(\theta),\alpha} : \ p \in \mathbf{G}_0\}$ which is given by $\mathbf{G}_0$ at this coordinate.
This can be done as $\cup\{ \eta_{p(\theta),\alpha} : \ p \in \mathbf{G}_0\} \in N_1$, since $\mathbf{G}_0, \overline{p}_0 \in N_1$. Then clearly $p'''_0 \leq p_{0,l}$ for each $l \in \omega$.

  Finally, for the last item of $\bref{p0}$ first recall that $\mathbb{P}^*_{S,\overline{U}_*}$ is an $\omega_1$-closed dense subposet of $\mathbb{P}_{S,\overline{U}_*}$ by Remark $\ref{P*}$. Then if a countable increasing sequence in $\mathbb{P}^*_{S,\overline{U}_*}$ (where a first element stronger than $p'''_0$) decides more and more about the $\delta_{\bullet}$'th level of $\dot{T}$, then choosing $p'_0$ to be an upper bound will work (e.g. choose an enumeration $\langle \dot{t}_i: \ i \in \omega \}$ of the $\delta_\bullet$'th level of $\dot{T}$, let $\langle s_i: \ i \in \omega \rangle$ enumerate $\dot{T}_{<\delta_{\bullet}}$ in type $\omega$, and let $r_j$ decide whether the $j$'th ordered pair in the countable set $\{ s_i: \ i \in \omega \} \times \{ \dot{t}_i: i \in \omega \}$ is in $\leq_{\dot{T}}$).

The next step is to construct the $\overline{p}_i$'s ($i>0$) and the $q_n$'s. This will be done simultaneously by induction. The induction is carried out in $V$, but each step can be done in $N_1$, which will guarantee that each $\overline{p}_n \in N_1$. 

It is straightforward to check that choosing $q_0 = p_*$ would satisfy our requirements, as e.g. $p_{0,0} = p_* \up \overline{U}_*$.
Then fixing $n>0$, and assuming that $\overline{p}_{i}, q_i$ are constructed for each $i<n$, first we construct $q_n$. 
Recall that $q_{n-1} \up (\overline{U}_* + \overline{W}_{n-1}) \in \mathbf{G}_{n-1}$ (by $\bref{qv} / \ref{qgen}$). 

Recall the definition of the set $E$ $\eqref{eDf}$, and let 
\[ E_{n-1} = \{ e \in E: \ \forall i<n \  \langle \varrho_i,\xi_i \rangle \notin e\}. \]
Using that $p_* \in \mathbb{P}_{S, \overline{U}_* + \overline{\id}_{S'}}$ forced that $\dot{b}$ is not an element of $V[G_{S, \overline{U}_* + \overline{W}_{n-1}}]$, i.e. there is no $\mathbb{P}_{S, \overline{U}_* + \overline{W}_{n-1}}$-name of it, we argue that
 \[ \begin{array}{rl}
 D =\{ &  p \in \mathbb{P}_{S, \overline{U}_* + \overline{W}_{n-1}}: \ \exists e, e' \in E_{n-1}
 \ (p \tieconcat e \leq q_{n-1}, \ p \tieconcat e' \leq q_{n-1}) \wedge \\
& (\exists \delta< \omega_1, \ t \neq t' \in \dot{T}_{\leq \delta} \sm \dot{T}_{<\delta}: \ \ (p \tieconcat e \Vdash t \in \dot{b}) \wedge  (p \tieconcat e' \Vdash t' \in \dot{b})) 
    \}
\end{array} \]
  is dense in $\mathbb{P}_{S, \overline{U}_* + \overline{W}_{n-1}}$ under $q_{n-1} \up (\overline{U}_* + \overline{W}_{n-1})$. Indeed, assume on the contrary that $q' \in \mathbb{P}_{S, \overline{U}_* + \overline{W}_{n-1}}$, $q' \leq q_{n-1} \up (\overline{U}_* + \overline{W}_{n-1})$ is such that that $D$ has no element under $q'$. Now for every $\delta< \omega_1$, consider the set
  \[ \begin{array}{rl}
  D_\delta =\{ & p\in \mathbb{P}_{S,\overline{U}_* + \overline{W}_{n-1}}: \ \ (p \leq q') \ \wedge \ (\exists e \in E_{n-1}: \ [ p \tieconcat e \leq q_{n-1}] \wedge  \\ 
  & \wedge \ [\exists t_{p,e,\delta} \in \dot{T}_{\leq \delta} \sm \dot{T}_{<\delta}: \ p \tieconcat e \ \Vdash t_{p,e,\delta} \in \dot{b} ]) \} 
  \end{array},\]
which is dense under $q'$ in $ \mathbb{P}_{S,\overline{U}_* + \overline{\id}{S'}}$. Now since for each $\delta < \omega_1$ the sets $D$ and $D_\delta$ are disjoint, for $p \in D_\delta$ the witnessing $t_{p,e,\delta}$ doesn't depend on $e$, therefore $q' \wedge q_{n-1}$ 
forces that $\dot{b}$ is in $V[\mathbf{G}_{S, \overline{U}_* + \overline{W}_{n-1}}]$ (i.e. forces that the $\mathbb{P}_{S, \overline{U}_* + \overline{W}_{n-1}}$-name $\{\langle p, t_{p,\delta} \rangle : \ \ p \in D_\delta, \ \delta < \omega_1 \}$ and $\dot{b}$ are equal). \\ Then as our set $D \in N_0$ is indeed dense we have that there exists a condition $q'' \in \mathbf{G}_{n-1} \cap D$, witnessed by $t \neq t'$ and $e,e'$. Finally, if $t \in B_n$ then define $q_{n} = q'' \tieconcat e'$, otherwise we can let $q_n = q'' \tieconcat e$, which are both stronger conditions than $q_{n-1}$ by the definition of $D$. It is straightforward to check $\bref{qv}$.

As $q_n$ is already defined (and so are $\overline{p}_{i}, q_i$ for each $i<n$), we turn to the definition of $\overline{p}_n$, which we will do similarly to that of $\overline{p}_0$. 
Let $p_{n,0} = q_n \upharpoonright (\overline{U}_* + \overline{w}_{n}$), assume that $p_{n,0}, p_{n,1}, \dots, p_{n,l-1}$ are already chosen. \\ If $\varepsilon_{l-1} \notin E_{n-1}$,
then $p_{n,l} = p_{n,l-1}$, otherwise proceed as follows.
Choose the sequence $\overline{e} = \overline{e}(n,l-1) = \langle e_i: \ 1\leq i \leq n+1 \rangle \in E^{n+1 \sm \{0 \}}$ and the sequence $\overline{m} = \overline{m}(n,l-1) = \langle m_i: \ i \leq n \rangle \in \omega^{n+1}$ with the property 
\begin{enumerate}[label = \arabic*)]
	\item \label{ek} $e_{n+1} = \varepsilon_{l-1}$ and $m_{n} = l-1$, 
	\item \label{mk} for each $i < n+1$
	\beeq \label{ekmk} J_{m_i} = D \ \  \wedge \ \ \text{"}e_i = \ (e_{i+1}  \ \text{ plus } \ (\eta_{p_{i,m_i}(\varrho_i),\xi_i} \text { attained on } \langle \varrho_i, \xi_i \rangle))\text{"}.  \eeq 
\end{enumerate}
Provided that the $e_j$'s are defined for $j>i$, and as well each $m_j$ for $j \geq i$, let $e_i \in E$ be the element with $u_{e_i} = u_{e_{i+1}} \cup \{ \langle \varrho_{i}, \xi_{i} \rangle \}$, $\overline{\eta}_{e_i} \supseteq \overline{\eta}_{e_{i+1}}$, $\eta_{e,\varrho_{i}, \xi_{i}} = \eta_{p_{i,m_{i}}(\varrho_i), \xi_i}$, and let $m_{i-1} = k(D,e_{i})$. Observe that by our procedure, and by the definition of the function $k$ $\eqref{kdf}$ we have $e_1 = \varepsilon_{m_0}$, and also 
\beeq \label{e1} \eta_{e_1,\varrho_n,\xi_n} = \eta_{p_{n,l-1}(\varrho_n),\xi_n}. \eeq
At some point later we will use the following fact, hence it is worth to note that for each $i$, $1 \leq i \leq n$
\beeq \label{kezdosz} \overline{e}(i,m_i) \subseteq  \overline{e}(n,l-1), \text{ and } \overline{m}(i,m_i) \subseteq \overline{m}(n,l-1).   \eeq
Finally consider the condition $r_{m_0}$ (from $\circledast_{\ref{p'1}}- \circledast_{\ref{p'l}}$): if $r_{m_0}\tieconcat e_1$ is a not a condition in $\mathbb{P}_{S,\overline{U}_* + \id \up S'}$, then let $p_{n,l} = p_{n,l-1}$, otherwise first define the auxiliary condition \beeq r_\bullet= g(r_{m_0}\tieconcat e_1, D), \eeq
  and note that in this case $\eta_{(r_{m_0}\tieconcat e_1)(\varrho_n),\xi_n} = \eta_{p_{n,l-1}(\varrho_n),\xi_n}$ by $\eqref{e1}$, therefore by the properties of $g$ we obtain
 \beeq \label{xinkiterj} \eta_{r_\bullet(\varrho_n),\xi_n} \supseteq \eta_{p_{n,l-1}(\varrho_n),\xi_n}. \eeq  
Recall that $p_{n,l-1}\up \overline{U}_* \in \mathbf{G}_0$ by our induction hypotheses $\bref{pvonas}$, and  it can be seen from the construction of $p_{0,j}$'s that in this case $p_{0,m_0+1} = r_\bullet \up \overline{U}_* \in \mathbf{G}_0$. Therefore by $\eqref{xinkiterj}$ we have that $(r_\bullet \up \overline{U}_* + \overline{w}_n) \wedge p_{n,l-1} $ is a condition in $\mathbb{P}_{\overline{U}_* + \overline{w}_n}$, and
let 
\[ p_{n,l} = (r_\bullet \up \overline{U}_* + \overline{w}_n) \wedge p_{n,l-1}. \]

 Then clearly $p_{n,l} \leq p_{n,l-1}$, and $p_{n,l} \up \overline{U}_* \in \mathbf{G}_0$.
 From $\bref{pvonas}$ it only remained to check that $\ref{N1ben}$ and $\ref{pngen}$ also hold. Since the whole construction of $\overline{p}_n$ took place in $N_1$ ($k \in N_1$ and so is the enumeration $\langle \langle J_i, \varepsilon_i \rangle: \ i \in \omega \rangle$, $g \in N_0$), $\overline{p}_n \in N_1$ obviously follows.
 Verifying the genericity of $\mathbf{G}_n$ goes similarly as of $\mathbf{G}_0$. Let $D \subseteq \mathbb{P}_{S, \overline{U}_* + \overline{W}_n}$, $D \in N_0$ be a fixed dense set, and $e' \in E$ be the empty sequence. Now, if we choose $l$  so that $J_{l-1} =  D' = \{ p \in \mathbb{P}_{S, \overline{U}_* + \overline{\id}_{S'}}: \ p \up \overline{U}_* + \overline{W}_n \in D\}$, $\varepsilon_{l-1} = e'$, then it follows from the construction of $p_{k,j}$'s, that of  $\overline{m} = \overline{m}(n,l-1)$ and $\overline{e} = \overline{e}(n,l-1)$, and from $\eqref{kezdosz}$ that 
 \[ p_{i,m_{i}+1} = (r_\bullet \up \overline{U}_* + \overline{w}_i) \wedge p_{i,m_i}  \ \text{ if }1 \leq i \leq n\text{,} \]
 and
 \[ p_{0,m_0+1} = g(r_{m_0} \tieconcat e_1) \up \overline{U}_*, \]   
 therefore
 \[ \bigwedge_{i \leq n} p_{i,m_i} \leq g(r_{m_0} \tieconcat e_1) \up (\overline{U}_* + \overline{W_n}) \in D'. \]



\ep(Claim $\ref{focl}$)

\ep (Lemma $\ref{fol}$)

\bl \label{lemma}
Let $T \in V[\mathbf{G}_{S,\overline{U}_*}]$ be a Kurepa tree, $S' \subseteq S \cap S_\bullet$ ($S' \in V$), $\mathbf{G}^\circ_{\overline{\id}_{S'}-\overline{U}_*} \subseteq \mathbb{P}^\circ_{\overline{\id}_{S'}-\overline{U}_*}$ be generic over $V[\mathbf{G}_{S,\overline{U}_*}]$. 
Suppose that $b \in V[\mathbf{G}_{S,\overline{U}_*}][\mathbf{G}^\circ_{S',(\overline{\id}_{S'}-\overline{U}_*)}] \sm V[\mathbf{G}_{S,\overline{U}_*}]$ is a new branch of $T$, and suppose that $\gamma \geq \kappa$ is a cardinal, and for each $\theta \in S'$ the inequality $|\theta \sm U^*_\theta| \geq \gamma$ holds. Then the filter $\mathbf{G}^\circ_{\overline{\id}_{S'}-\overline{U}_*}$ adds at least $|\gamma|$-many new branches to $T$.
\el
\bp
W.l.o.g. we can assume that $T \subseteq \omega_1$, and $\lambda$ is a cardinal (in $V[\mathbf{G}_{S,\overline{U}_*}]$). First we will choose a system $\overline{W}_0 = \langle W_{0,\theta}: \ \theta \in S' \rangle \in \prod_{\theta \in S'} \mathcal{P}(\theta)$ with  ($\forall \theta \in S'$) $|W_{0,\theta}| < \kappa$, and $b \in V[\mathbf{G}_{S,\overline{U}_*}][\mathbf{G}^\circ_{\overline{W}_0}]$:
 since $b \in V[\mathbf{G}_{S,\overline{U}_*}][\mathbf{G}^\circ_{\overline{\id}_{S'}-\overline{U}_*}]$, $S' \in V$ we can use Lemma $\ref{mindegyforsz}$  and obtain  that $b \in V[\mathbf{G}_{S,\overline{U}_*}][\mathbf{G}^\circ_{\overline{\id}_{S'}-\overline{U}_*}] = V[\mathbf{G}_{S, \overline{U}_* + \overline{\id}_{S'}}]$. And because $b \subseteq \mathcal{H}(\omega_1)^V$, applying Lemma $\ref{her}$ with $S$, and $\overline{U} = \overline{U}_* + \overline{\id}_{ S'}$, there exists $S_* \subseteq S$, $\overline{W}_* \in \prod_{S_* \sm S'} \mathcal{P}(U_\theta) \times \prod_{\theta \in S_* \cap S'} \mathcal{P}(\theta)$ with 
 \[ b \in V[\mathbf{G}_{S_*, \overline{W}_*}] \subseteq V[\mathbf{G}_{S, \overline{U}_* + \overline{W}_*}] = V[\mathbf{G}_{S, \overline{U}_*}][\mathbf{G}^\circ_{\overline{W}_* - \overline{U}_*}] , \]
  where  $|S_*| < \kappa$, and $|W^*_\theta| < \kappa$ for each $\theta \in S_*$. Then fixing $\overline{W}_0 \in \prod_{\theta \in S'} \mathcal{P}(\theta)$  so that $W_{0,_\theta} = W^*_{\theta} \sm U^*_\theta$ if $\theta \in S_*$, and $W_{0,\theta} = \emptyset$ for $\theta \in S \sm S_*$ has the required  properties.
   
   Now, as $|W_{0,\theta}| < \kappa \leq \gamma$, and $\gamma \leq |\theta \sm U^*_\theta|$ for each $\theta \in S'$ we can fix for each $\alpha < \gamma$ a system $\overline{W}_\alpha = \langle W_{\alpha, \theta}: \ \theta \in S' \rangle \in \prod_{\theta \in S'} \mathcal{P} ( \theta \sm U^*_\theta)$ such that for every $\theta \in S'$
   \begin{enumerate}[label = (\roman*)]
   	   	\item $W_{\alpha, \theta} \cap W_{\beta, \theta} = \emptyset$ for every $\alpha< \beta < \gamma$,
   	   	\item $|W_{0,\theta}| = |W_{\alpha,\theta}|$ for each $\alpha < \gamma$.
   \end{enumerate}
   For each $0<\alpha < \gamma$ define the bijections 
   \[ \pi_\alpha : \bigcup_{\theta \in S'} \{\theta\} \times W_{0,\theta} \to \bigcup_{\theta \in S'} \{\theta\} \times W_{\alpha,\theta} \]
   where $\pi_\alpha \up \{\theta\} \times W_{0,\theta}$ is a bijection to $\{\theta\} \times W_{\alpha,\theta}$. Then clearly each $\pi_{\alpha}$ induces an automorphism $\hat{\pi}_{\alpha} \in V[\mathbf{G}_{S, \overline{U}_*}]$ of $\mathbb{P}^\circ_{\overline{W}_0}$ and $\mathbb{P}^\circ_{\overline{W}_\alpha}$. Moreover, $\hat{\pi}_\alpha$ induces a natural operation $\hat{\pi}_{\alpha}^*$ from the class of  $\mathbb{P}^\circ_{\overline{W}_0}$-names to the class of $\mathbb{P}^\circ_{\overline{W}_\alpha}$-names.
   Now fix a $\mathbb{P}^\circ_{\overline{W}_0}$-name $\dot{b}_0 \in V[\mathbf{G}_{S, \overline{U}_*}]$ for our new branch $b \in  V[\mathbf{G}_{S, \overline{U}_*}][\mathbf{G}^\circ_{\overline{W}_0}]$, and choose an element $p_\bullet \in \mathbb{P}^\circ_{\overline{W}_0}$ forcing that $\dot{b}_0$ is a new branch, i.e.
   \beeq \label{pbul} V[\mathbf{G}_{S, \overline{U}_*}] \models \ \ \ p_\bullet \Vdash \ \dot{b}_0 \in \mathcal{B}(T) \sm \mathcal{B}^{V[\mathbf{G}_{S, \overline{U}_*}]} (T). \eeq 
   Let $\mathbb{P}^\circ_\bullet = \mathbb{P}^\circ_{\sum_{\alpha < \gamma} \overline{W}_\alpha}$, i.e. adding the branches $\bigcup_{\alpha \in \gamma} W_{\alpha,\theta}$ to $\name{T}_\theta$ for each $\theta \in S'$, which is of course equal to the countably supported product of $\mathbb{P}^\circ_{\overline{W}_\alpha}$'s ($\alpha < \gamma$), and let $\mathbf{G}^\circ_\bullet$ denote  the generic filter $\mathbf{G}^\circ_{\overline{\id}_{S'} -  \overline{U}_*} \cap \mathbb{P}^\circ_\bullet$.
   
   We will show that in $V[\mathbf{G}_{S, \overline{U}_*}][\mathbf{G}^\circ_\bullet] \subseteq V[\mathbf{G}_{S, \overline{U}_*}][\mathbf{G}^\circ_{\overline{\id}_{S'}-\overline{U}_*}] $ there are at least $\gamma$-many new branches of $T$, i.e.
   \[  \left| \mathcal{B}(T) \cap \left(V[\mathbf{G}_{S, \overline{U}_*}][\mathbf{G}^\circ_\bullet] \sm V[\mathbf{G}_{S, \overline{U}_*}] \right) \right| \geq \lambda, \] 
   by arguing that 
   \begin{enumerate}[label = $\otimes_{\arabic*}$, ref =\arabic*]
   	\item \label{paforsz} for any $\alpha < \gamma$  (in $V[\mathbf{G}_{S,\overline{U}_*}]$)
   	\[ \hat{\pi}_{\alpha}(p_\bullet) \Vdash_{\mathbb{P}_\bullet^\circ} \hat{\pi}^*_{\alpha}(\dot{b}_0) \notin V[\mathbf{G}_{S,\overline{U}_*}][\mathbf{G}^\circ_{\bullet, <\alpha}] \] (where $\mathbf{G}^\circ_{\bullet, <\alpha}$ stands for $\mathbf{G}^\circ_{\bullet} \cap \mathbb{P}^\circ_{\sum_{\beta < \alpha} \overline{W}_\beta}$), and 
   	 \item \label{kett} $|\{ \alpha < \gamma: \ \hat{\pi}_\alpha(p_\bullet) \in \mathbf{G}_\bullet^\circ  \}| = \gamma $.
   \end{enumerate} 
This will complete the proof of Lemma $\ref{lemma}$.
 
 First we will prove $\otimes_{\ref{kett}}$, for which recall that we assumed that $\gamma$ is a cardinal, and choose a system of uncountable regular cardinals 
    $\{\rho_{\beta}: \ \beta < \chi < \gamma \}$, and a partition  $\langle I_\beta: \ \beta < \chi \rangle$ of $\gamma$ with $\otp(I_\beta) = \rho_{\beta}$ for each $\beta <\chi$ (i.e. $I_\beta \cap I_{\delta} = \emptyset$ for $\beta<\delta <\rho$, and $\bigcup_{\beta < \rho} I_\beta = \gamma$).
    Then it is enough to verify
\beeq \label{Ibe} (\forall \beta < \chi) \ \ \ |\{ \alpha \in I_\beta : \ \hat{\pi}_\alpha(p_\bullet) \in \mathbf{G}_\bullet^\circ  \}| = \rho_\beta, \eeq
which can be seen by a standard density argument: Fix $\beta < \varrho$, $\alpha \in I_\beta$, then it suffices to show that 
\[ D_{\beta,\alpha} = \{ p \in \mathbb{P}^\circ_\bullet: \ p \leq \hat{\pi}_\delta(p_\bullet) \ \text{ for some } \ \delta> \alpha, \ \delta \in I_\beta \} \text{ is dense,}\]
which obviously holds by the regularity of the uncountable $\rho_\beta = |I_\beta|$ (since for $\delta \in I_\beta$ we have $\hat{\pi}_\delta(p_\bullet) \in \mathbb{P}^\circ_{\overline{W}_\delta}$,  $\mathbb{P}^\circ_\bullet$ is the countably supported product of $\mathbb{P}^\circ_{\overline{W}_\alpha}$'s ($\alpha < \gamma$), and $I_\beta \subseteq \gamma$).

For $\otimes_{\ref{paforsz}}$ first consider $\mathbb{P}_\bullet^\circ$ as the product of $\mathbb{P}^\circ_{\sum_{\beta < \gamma, \beta \neq \alpha} \overline{W}_\beta}$ and $\mathbb{P}^\circ_{\overline{W}_\alpha}$. We will need the following claim.
\bc \label{agcl} For each $p \in \mathbb{P}^\circ_{\overline{W}_\alpha}$, $p \leq \hat{\pi}_{\alpha}(p_\bullet)$ there exist $q_0, q_1 \in  \mathbb{P}^\circ_{\overline{W}_\alpha}$ $q_0, q_1 \leq p$, and the  incomparable elements $t_0, t_1$ of the tree $T$ such that 
\[ V[\mathbf{G}_{S,\overline{U}_*}][\mathbf{G}^\circ_{\bullet, \gamma \setminus \{\alpha \} }] \models \ \ \ (q_i \Vdash_{\mathbb{P}_{\overline{W}_\alpha}^\circ} t_i \in \hat{\pi}^*_{\alpha}(\dot{b}_0)) \ \  \text{ for each } i \in \{0,1\},  \]
where $\mathbf{G}^\circ_{\bullet, \gamma \setminus \{\alpha \} } = \mathbf{G}^\circ_\bullet  \cap \mathbb{P}^\circ_{\sum_{\beta < \gamma, \beta \neq \alpha} \overline{W}_\beta}$.
\ec
Before proving the claim we verify that $\otimes_{\ref{paforsz}}$ follows from it. In fact 
\[ \hat{\pi}_{\alpha}(p_\bullet) \Vdash_{\mathbb{P}_\bullet^\circ} \ \hat{\pi}^*_{\alpha}(\dot{b}_0) \notin V[\mathbf{G}_{S,\overline{U}_*}][\mathbf{G}^\circ_{\bullet, \gamma \setminus \{\alpha \} }]. \]
Since $\mathbf{G}^\circ_{\bullet} \subseteq \mathbb{P}^\circ_\bullet$ is generic over $V[\mathbf{G}_{S,\overline{U}_*}]$, and $\mathbb{P}^\circ_\bullet$ can be identified with 
\[ \left(\mathbb{P}^\circ_{\sum_{\beta < \gamma, \beta \neq \alpha} \overline{W}_\beta} \right) \times \mathbb{P}^\circ_{\overline{W}_\alpha},\]
by \cite[Lemma V.1.1]{Ku2013} $\mathbf{G}^\circ_{\bullet, \gamma \setminus \{\alpha \} } = \mathbf{G}^\circ_\bullet  \cap \mathbb{P}^\circ_{\sum_{\beta < \gamma, \beta \neq \alpha} \overline{W}_\beta}$ is generic over $V[\mathbf{G}_{S,\overline{U}_*}]$, and $\mathbf{G}^\circ_{\bullet, \alpha } =  \mathbf{G}^\circ_\bullet  \cap \mathbb{P}^\circ_{\overline{W}_\alpha}$ is generic over $V[\mathbf{G}_{S,\overline{U}_*}][\mathbf{G}^\circ_{\bullet, \gamma \setminus \{\alpha \} }]$. For each branch  $c \in V[\mathbf{G}_{S,\overline{U}_*}][\mathbf{G}^\circ_{\bullet, \gamma \setminus \{\alpha \} }]$ of $T$ define (in $V[\mathbf{G}_{S,\overline{U}_*}][\mathbf{G}^\circ_{\bullet, \gamma \setminus \{\alpha \} }]$)
\[ D_c = \{ q \in \mathbb{P}^\circ_{\overline{W}_\alpha}: \ \exists t \in T \sm c \ \text{ such that }  q \Vdash_{\mathbb{P}_{\overline{W}_\alpha}^\circ}  t \in \hat{\pi}^*_{\alpha}(\dot{b}_0) \} 
, \]
which is dense under $\hat{\pi}_{\alpha}(p_\bullet)$ by Claim $\ref{agcl}$, since for a fixed $p \in \mathbb{P}^\circ_{\overline{W}_\alpha}$ at most one $t_i$ can be in the branch $c$.

\bp(Claim $\ref{agcl}$)
 First we argue that the statement holds in $V[\mathbf{G}_{S,\overline{U}_*}]$, i.e. for each $p \in \mathbb{P}^\circ_{\overline{W}_\alpha}$, $p \leq \hat{\pi}_{\alpha}(p_\bullet)$ there exist $q_0, q_1 \in  \mathbb{P}^\circ_{\overline{W}_\alpha}$ $q_0, q_1 \leq p$, and the  incomparable elements $t_0, t_1$ of the tree $T$ such that 
 \beeq \label{qii} V[\mathbf{G}_{S,\overline{U}_*}] \models \ \ \ (q_i \Vdash_{\mathbb{P}_{\overline{W}_\alpha}^\circ} t_i \in \hat{\pi}^*_{\alpha}(\dot{b}_0)) \ \  \text{ for each } i \in \{0,1\}.  \eeq
 Now $\eqref{pbul}$ implies that 
 \[ V[\mathbf{G}_{S, \overline{U}_*}] \models \ \ \ \hat{\pi}_\alpha(p_\bullet) \Vdash_{\mathbb{P}^\circ_{\overline{W}_\alpha}} \ \hat{\pi}_\alpha^*(\dot{b}_0) \in \left(\mathcal{B}(T) \sm  \mathcal{B}^{V[\mathbf{G}_{S, \overline{U}_*}]} (T)\right) \]
 since $\dot{b}_0 \in V[\mathbf{G}_{S, \overline{U}_*}]$ is a $\mathbb{P}^\circ_{\overline{W}_0}$-name and $T \in V[\mathbf{G}_{S, \overline{U}_*}]$. Suppose that $p \leq \hat{\pi}_\alpha(p_\bullet)$ is a counterexample, but then for the set
 \[ b' = \{ t \in T: \ \exists q \in \mathbb{P}^\circ_{\overline{W}_\alpha}, \ q \leq p \ \text{ s.t. } q \Vdash t \in \hat{\pi}_\alpha^*(\dot{b}_0)  \} \in V[\mathbf{G}_{S, \overline{U}_*}] \]
 we have $p \Vdash \hat{\pi}_\alpha^*(\dot{b}_0) = b'$ (since $\hat{\pi}_\alpha(p_\bullet)$ forced that $\hat{\pi}_\alpha^*(\dot{b}_0)$ is a cofinal branch in $T$), a contradiction.
 Finally, fixing $p \leq \hat{\pi}_{\alpha}(p_\bullet)$, if $q_0, q_1 \in  \mathbb{P}^\circ_{\overline{W}_\alpha}$ $q_0, q_1 \leq p$, and the  incomparable elements $t_0, t_1 \in T$ are such that $\eqref{qii}$ holds, then
 \[ V[\mathbf{G}_{S,\overline{U}_*}][\mathbf{G}^\circ_{\bullet, \gamma \setminus \{\alpha \} }] \models \ \ \ (q_i \Vdash_{\mathbb{P}_{\overline{W}_\alpha}^\circ} t_i \in \hat{\pi}^*_{\alpha}(\dot{b}_0)) \ \  \text{ for each } i \in \{0,1\},  \]
 since if $q_i \in \mathbf{H} \subseteq \mathbb{P}_{\overline{W}_\alpha}^\circ$ is generic over $V[\mathbf{G}_{S,\overline{U}_*}][\mathbf{G}^\circ_{\bullet, \gamma \setminus \{\alpha \} }]$, and $t_i \notin \hat{\pi}^*_{\alpha}(\dot{b}_0)[\mathbf{H}]$ (for some $i \in \{0,1\}$), then $\mathbf{H}$ is generic over $V[\mathbf{G}_{S,\overline{U}_*}]$ too, and the same holds in $V[\mathbf{G}_{S,\overline{U}_*}][\mathbf{H}]$.
\ep

It is left to argue why Lemma $\ref{fol}$ and Lemma $\ref{lemma}$ complete the proof of Theorem $\ref{fot}$ (and Theorem $\ref{fot2}$). Suppose that $T \in V[\mathbf{G}]$ is a Kurepa tree (where $\mathbf{G} \subseteq \mathbb{P} = \mathbb{P}_{S_\bullet^+, \overline{\id}_{S_\bullet^+}}$ is generic), and assume on the contrary that $|\mathcal{B}^{V[\mathbf{G}]}(T)| \notin S_\bullet$. We can also assume that $T \subseteq \mathcal{H}(\omega_1)^V$, and by Lemma $\ref{her}$ there exists $S_* \subseteq S^+_\bullet$, $|S_*| < \kappa$,  $\overline{W}_* = \langle W^*_\theta: \ \theta \in S_* \rangle \in \prod_{\theta \in S_*} [\theta]^{<\kappa}$ such that $T \in V[\mathbf{G}_{S_*, \overline{W}_*}]$. For estimating $(2^{\omega_1})^{V[\mathbf{G}_{S_*, \overline{W}_*}]}$ first a straightforward calculation yields that
$|\mathbb{P}_{S_*, \overline{W}_*}| < \kappa$: Since $|\mathbb{P}_{S_*,\langle \emptyset: \ \theta \in S_* \rangle}| = (|S_*||\omega_1|)^\omega$ which is either $(\omega_1 \cdot \omega_1)^\omega = \omega_1 < \omega_2$ (if $\kappa = \omega_2$, by $\mathbf{CH}$), or $\gamma^\omega < \kappa$ (for some $\gamma < \kappa$, if $\kappa$ is inaccessible). Thus recalling the definition of $\mathbb{Q}_{\theta, W^*_\theta}$'s, the fact $\sum_{\theta \in S_*} |W^*_\theta| < \kappa$ as $\kappa$ is regular, and $\sup W^*_\kappa < \kappa$ (if $\kappa \in S_*$) we have the following (in both cases regardless of whether $\kappa = (\omega_2)^V$, or an inaccessible)
\[ |\mathbb{P}_{S_*, \overline{W}_*}| = |\mathbb{P}_{S_*, \langle \emptyset: \ \theta \in S_* \rangle}| \cdot \left((\omega_1) \cdot \left(\sum_{\theta \in S_* \sm \{ \kappa \}} |W^*_\theta|\right)\right)^{\omega} \cdot (|W^*_\kappa| \cdot \sup W^*_\kappa)^\omega < \kappa. \]
At this point we have to discuss the two cases (i.e. whether $\kappa \in S_\bullet$) differently, arguing that in both cases there are branches outside $V[\mathbf{G}_{S_*, \overline{W}_*}]$. 

If $\kappa = \omega_2 \in S_\bullet$, then as 
\[ V \models \ |\mathbb{P}_{S_*, \overline{W}_*}|^{\omega_1 \cdot  |\mathbb{P}_{S_*, \overline{W}_*}|} = \omega_2, \]
we have
\[ V[\mathbf{G}_{S_*, \overline{W}_*}] \models \ 2^{\omega_1} = \omega_2,\]
therefore as $|\mathcal{B}^{V[\mathbf{G}]}(T)| \notin S_\bullet$, there are branches of $T$ in $V[\mathbf{G}]$ not in $V[\mathbf{G}_{S_*, \overline{W}_*}]$. On the other hand, if $\kappa \notin S_\bullet$ is inaccessible, then we obtain that 
\[ V[\mathbf{G}_{S_*, \overline{W}_*}] \models \ |\mathcal{B}(T)| \leq 2^{\omega_1} < \kappa, \]
and as $\kappa$ remains a cardinal in $V[\mathbf{G}]$ (by Claim $\ref{card}$), and 
\[ V[\mathbf{G}] \models \ |\mathcal{B}(T) \cap V[\mathbf{G}_{S_*, \overline{W}_*}] | = \omega_1, \]
we conclude that this case there also must be branches of $T$ not in $V[\mathbf{G}_{S_*, \overline{W}_*}]$ as $T$ is a Kurepa tree in $V[\mathbf{G}]$.
Now let $\overline{R} \in \prod_{\theta \in S_\bullet^+ \sm S_\bullet} \mathcal{P}(\theta)$, $R_\theta = \theta \sm W^*_\theta$, then 
\[ \mathbb{P} = \mathbb{P}_{S_\bullet^+, \overline{\id}_{S_\bullet^+}} \simeq (\mathbb{P}_{S_*, \overline{\id}_{S*}- \overline{R}}) \times (\mathbb{P}_{S_* \cap (S_\bullet^+ \sm S_\bullet), \overline{R}}) \times (\mathbb{P}_{S^+_\bullet \sm S_*, \overline{\id}_{S_\bullet^+ \sm S_*}}), \]
and there are no new sequences of type $\omega$ in $V[\mathbf{G}]$ (by Claim $\ref{ome}$), and the second component is $\omega_1$-closed, the third component has an $\omega_1$-closed dense subset (which thus remain $\omega_1$-closed in $V[\mathbf{G}_{S_*, \overline{\id}_{S_*} - \overline{R}}]$) we obtain that each branch of $T$ is added by $\mathbf{G}_{S_*, \overline{\id}_{S_*} - \overline{R}} = \mathbf{G} \cap \mathbb{P}_{S_*, \overline{\id}_{S*} - \overline{R}}$ (since an $\omega_1$-closed forcing do not add new branches to Kurepa trees \cite[Lemma V.2.26]{Ku2013}). We only have to derive a contradiction from
\[ V[\mathbf{G}_{S_*, \overline{\id}_{S_*} - \overline{R}}] \models \ |\mathcal{B}(T)| \notin S_\bullet. \]
Now  letting $\partial=|\mathcal{B}^{V[\mathbf{G}_{S_*, \overline{\id}_{S_*}- \overline{R}}]}(T)|  \notin S_\bullet$,  $S_{*}^- = S_* \cap S_\bullet  \cap \partial$, $S_*^+ = (S_* \cap S_\bullet) \sm S_*^-$
by Lemma $\ref{mindegyforsz}$ we have 
\[ V[\mathbf{G}_{S_*, \overline{\id}_{S_*}- \overline{R}}] = V[\mathbf{G}_{S_*,\overline{W}_* + \overline{\id}_{S_*^-}}][\mathbf{G}^{\circ}_{\overline{\id}_{S_*^+}- \overline{W}_*}].\]



As $\partial \notin S_*^-, S_*^+$, it is enough to prove that in $V[\mathbf{G}_{S_*,\overline{W}_* + \overline{\id}_{S_*^-}}]$ there are less than $\partial$-many branches of $T$, because if $\mathbf{G}^{\circ}_{\overline{\id}_{S_*^+}- \overline{W}_*}$ adds new branches, then adds $\min(S_*^+)$-many new branches by Lemma $\ref{lemma}$ (since each $|W^*_\theta|<\kappa \leq \min(S_\bullet) \leq \min(S_*^+)$). 

Now if $\partial = \kappa$, then $S_*^- = \emptyset$, we are done, so we can assume that $\partial > \kappa$, and $\sup S_*^- \geq \kappa$. As $|S_*| < \kappa$ (in $V$), and our conditions ($\underline{\text{Case1}}/\ref{harom}$, or $\underline{\text{Case2}}/\ref{harom'}$) states that then $\sup(S_* \cap S_\bullet \cap \partial) \in S_\bullet$ implying $\sup S_*^- < \partial$. Therefore using that $W^*_\theta \subseteq \theta$  we get $\sum_{\theta \in S_*^-} |W^*_\theta| \leq |\sup S_*^-|^2 < \partial$. Now by Lemma $\ref{fol}$ for each branch $b$ of $T$ 
in $V[\mathbf{G}_{S_*,\overline{W}_* + \overline{\id}_{S_*^-}}] = V[\mathbf{G}_{S_*,\overline{W}_*}][\mathbf{G}^\circ_{(\overline{\id}_{S_*^-}) - \overline{W}_*}]$ there exist $\theta_0, \theta_1, \dots, \theta_{n-1}$, $U^\bullet_{\theta_0}, U^\bullet_{\theta_1}, \dots, U^\bullet_{\theta_{n-1}}$ finite such that  $b \in V[\mathbf{G}_{S_*,\overline{W}_*}][\mathbf{G}^\circ_{\overline{U}_\bullet}]$. Therefore, as $|\mathbb{P}^\circ_{\overline{U}_\bullet}| = \omega_1^n = \omega_1$, counting the nice $\mathbb{P}^\circ_{\overline{U}_\bullet}$-names of subsets $T$ for each possible $n$, sequence of $\theta$'s, and $\overline{U}_\bullet$
\[ \mathcal{B}(T) \cap (V[\mathbf{G}_{S_*,\overline{W}_*}][\mathbf{G}^\circ_{(\overline{\id}_{S_*^-}) - \overline{W}_*}] \sm V[\mathbf{G}_{S_*,\overline{W}_*}]) \leq (|\sup S_*^-|^{<\omega} \cdot \omega_1^{\omega_1})^{V[\mathbf{G}_{S_*,\overline{W}_*}]} \leq \sup S_*^-, \]
which is smaller than $\partial$, a contradiction.

For $V[\mathbf{G}] \models \ 2^{\omega_1} = \lambda$ we only need to show that $2^{\omega_1} \leq \lambda$. But a similar straightforward calculation yields that $\mathbb{P} = \mathbb{P}_{S_\bullet^+, \overline{\id}_{S_\bullet^+}}$ is of cardinality $\lambda$, and then (using $\kappa$-cc and the equality $\lambda^{<\kappa} = \lambda$) by counting the possible nice names for subsets of $\omega_1$ we obtain the desired inequality.

\br If $S_\bullet$ also satisfies 
\beeq \label{nemelh} \forall \mu \in S_\bullet: \ \cf(\mu) < \kappa \ \to \ \mu^+ \in S_\bullet, \eeq
and $\mathbf{GCH}$ holds in $V$ then $S_\bullet \sm \{ \lambda \}$ is the spectrum for the Jech-Kunen trees in $V[\mathbf{G}]$. (A tree $T$ of height $\omega_1$ and power $\omega_1$ is a Jech-Kunen tree if $\omega_1 < |\mathcal{B}(T)| <2^{\omega_1}
$.) For more on Jech-Kunen trees see also \cite{JiSh:466}, \cite{JiSh:469}, \cite{JiSh:498}.
Note that $\mathbf{CH}$ in the final model implies that the product of countably many Jech-Kunen trees is a Jech-Kunen tree, so is the diagonal product of $\omega_1$-many Jech Kunen trees, hence $\eqref{nemelh}$ cannot be dropped.

One can obtain similar cardinal arithmetic conditions for $\Sp_\mu$ with $\mu$ large.
\er


 

\section{The necessity of the inaccessible cardinal}
In this section we prove that if $\omega_2$ is not an element of the spectrum, then $\omega_2$ is inaccessible in $L$. The idea of using transitive collapses of elementary submodels of constructible sets as nodes of a tree goes back to Solovay's original unpublished argument for the consistency strength of the negation of the Kurepa Hypothesis. Although the next proof is deemed to be well-known, for the sake of completeness we include the proof as there is probably no known source to cite.

\bt \label{szukst} Suppose that $\omega_2^V$ is a successor in $L$. Then there exists a Kurepa tree $T$ with $\mathcal{B}^V(T) = \omega_2$. 
\et
\bp
We will use an extension of $L$, an inner model between $L$ and $V$, what serves as the motivation for the following definition of relative constructibility, which can be found in e.g. \cite{kanamori2003higher}. 
\bd For a set $A$ define $L[A] = \bigcup_{\alpha \in ON} L_\alpha[A]$ by transfinite recursion as follows. $L_0[A] = \emptyset$, $L_{\alpha+1}[A] = \ddf_A(L_\alpha[A])$, and  $\alpha$ limit $L_\alpha[A] = \bigcup_{\beta < \alpha} L_\beta[A]$ (where $\ddf_Y(X)$ are the subsets of $X$ that can be defined in the structure $(X, \in \up (X \times X), Y \cap X)$ by parameters from $X$, see \cite[Chapter 1, §3]{kanamori2003higher}.
\ed
The following is standard easy exercise, but for the sake of completeness we include the proof.
\bc
There exists a set $A \subseteq \omega_1$ such that $\omega_1^{L[A]} = \omega_1$, $\omega_2^{L[A]} = \omega_2$.
\ec
\bp
If $\omega_2^V = (\lambda^+)^L$, where $|\lambda| = \omega_1$, then in a single subset $A$ of $\omega_1$ we can code a well-ordering of $\omega_1$ in type $\lambda$, and also for each $\alpha < \omega_1$ a well-ordering of $\omega$ in type $\alpha$ in the obvious fashion, and such that $L$ can read this coding (implying $\omega_1^{L[A]} = \omega_1$, $\omega_2^{L[A]} = \omega_2$): First let $\langle X_\alpha: \ \alpha \leq \omega_1 \rangle \in L$ be a set of pairwise disjoint sets of $\omega_1$ with $|X_\alpha|^L = \omega$ for each $\alpha < \omega_1$, and $|X_{\omega_1}|^L = \omega_1$, then for each $\alpha < \omega_1$ we can code the well ordering $X_\alpha$ in order type $\alpha$, and the well ordering of $X_{\omega_1}$ in type $\lambda$ in a subset $A'$ of $\bigcup_{\alpha \leq \omega_1} X_\alpha^2 \subseteq \omega_1^2$. Finally, taking the preimage of this set under a bijection $f \in L$ between $\omega_1$ and $\omega_1^2$, i.e. $A = f^{-1}(A')$ works. 
\ep

We have to recall a classical Lemma \cite[Theorem 3.3]{kanamori2003higher}. Recall that $\mathcal{L}_\in(R_A)$ stands for the (first-order) language of set theory extended by the unary predicate $R_A$.
\bl \label{conL} There is a sentence $\sigma \in \mathcal{L}_\in(R_A)$ such that for every transitive set $N$
\[ (N, \in, X \cap N)\models \sigma \text{ implies } N = L_\gamma[X] \text{ for some limit } \gamma. \]
In particular,
if $M \prec (L_\beta[X], \in, X \cap L_\beta[X])$, where $\beta$ is a limit ordinal and $\pi$ is the collapsing isomorphism from $M$ onto the transitive set $\ran(\pi)$, then the Mostowski collapse \[ \ran(\pi)= L_\gamma[\{ \pi(x): \ x \in M \cap X \}] \]
 for some $\gamma \leq \beta$.
\el
The following is immediate.
\bc \label{haszncor} For each infinite ordinal $\beta$ and $Y \subseteq L_\beta[X]$, if $Y \in L[X]$ and $X \subseteq L_\beta[X]$, then $\mu = (|\beta|^+)^{L[X]}$ implies $Y \in L_{\mu}[X]$. 
\ec
(Working in $L[X]$, if $Y \in L_\gamma[X]$, then let $M \prec L_\gamma[X]$ with $\{ Y \} \cup L_\beta[X] \subseteq M$, $|M| = |L_\beta[X]|$, and apply the lemma recalling that $\pi \up L_\beta[X]$ is the identity.) 

Now we can turn to the definition of the tree $T$, which will be defined by its branches.

Recall that there exists a definable well-order on $L[A]$, which is downward absolute to almost every initial segment of $L[A]$ (to the ones indexed by limit ordinals) \cite[Theorem 3.3]{kanamori2003higher}:
\bl \label{defrend} There exists a formula $\varphi \in \mathcal{L}_{\in}(R_A)$ (i.e. in the language of set theory extended with the unary relation symbol $A$) which define a well-ordering on $(L[A], \in, A)$, moreover if $\delta$ is a limit ordinal, $x,y \in L_\delta[A]$, then
\[ (L[A], \in, A) \models \varphi(x,y) \ \ \iff \ \ (L_\delta[A], \in, A \cap L_\delta[A]) \models \varphi(x,y).
 \]
\el
From now on '$x<_{L[A]}y$' abbreviates $\varphi(x,y)$.

We will take Skolem hulls many times, thus we need to introduce the following variant of this standard notion.
\bd \label{Skdef}Let $(M,\in, X, \partial)$, $M \subseteq L[A]$  be a set model of the language $\mathcal{L}_\in(R_A, c_\partial)$ with $\emptyset \in M$,  $M' \subseteq M$ such that the well-ordering formula $\varphi \in \mathcal{L}_\in(R_A)$ from Lemma $\ref{defrend}$ is absolute to $M$, i.e. 
\beeq \label{rab} (\forall x,y \in M): \ (L[A], \in, A) \models \varphi(x,y) \ \  \text{iff} \ \ (M, \in, X) \models \varphi(x,y), \eeq 
e.g. when $(M, \in, X) = (L_{\zeta}[A], \in, A \cap L_{\zeta}[A])$ for some limit ordinal $\zeta$. Then 
the Skolem-hull of $M'$ in $(M, \in, X, \partial)$ (in symbols, $\mathfrak{H}^{(M, \in, X, \partial)}(M')$)  is the closure of $M'$ under the functions $f^{(M, \in, X, \partial)}_{\psi}$ for each formula $\psi(v_0,v_1, \dots,v_{n_\psi}) \in \mathcal{L}_\in(R_A, c_\partial)$ with $n_\psi+1$ free variables, where the function $f^{(M, \in, X, \partial)}_{\psi}$ satisfies the following.
\[ f^{(M, \in, X, \partial)}_\psi: M^{n_\psi} \to M \] is defined so that for every $\langle x_1,x_2, \dots, x_{n_\psi} \rangle \in M^{n_\psi}$:
\[\text{ if } \exists y!\in M \ \text{ s.t. } (M,\in, X, \partial) \models \psi(y,x_1,x_2, \dots, x_{n_\psi}), \]
\[ \text{ then let } f^{(M, \in, X, \partial)}_\psi(x_1,x_2, \dots, x_{n_\psi})  \text{ be the unique such }y,\]  
\[ \text{otherwise let } f^{(M, \in, X, \partial)}_\psi(x_1,x_2, \dots, x_{n_\psi}) = \emptyset.\]
Then the fact that for each formula $\psi'$ we can define the formula saying that $y$ is the least $y$ (w.r.t. the well-order given by $\varphi$) satisfying $\psi'(y,x_1,x_2, \dots x_{n_{\psi'}})$ together with the  Tarski-Vaught criterion implies that the closure is an elementary submodel of $M$, in symbols, $M' \prec (M, \in, X, \partial)$.
\ed
Observe that this closure only depends on the isomorphism class of $(M, \in, X, \partial)$ by the absoluteness of the well-ordering formula $\varphi$ $\eqref{rab}$.

Choose $\xi < \omega_2$ such that
\beeq \label{xi} \xi \text{ is the minimal ordinal }(\forall \alpha < \omega_1)  \ \exists f_\alpha \in L_\xi[A]  \ \text{ bijection between } \omega \text{ and } \alpha \eeq
(which can be done due to Corollary $\ref{haszncor}$, in fact $\xi = \omega_1$, but we won't use this equality, hence we don't argue that).

Now we will define an operation which assigns for each $\delta \in [\xi, \omega_2)$ the ordinal $\delta' < \omega_2$ in the following way. We would like to choose $\delta'$ so that in $L_{\delta'}[A]$ it is true that for each set $x$ there exists a surjection from $\omega_1$ to $x$, and for $\delta'' \neq \delta'$ the structures $(L_{\delta'}[A], \in, A, \delta_0)$ and $(L_{\delta''}[A], \in, A, \delta_0)$ cannot be elementarily equivalent. 

\bd \label{delta'} Fix $\delta \in [\xi, \omega_2)$, and define $\delta'$ to be the least ordinal such that 
\begin{enumerate}[label=\alph*)]
	\item \label{d'1} $\delta \in L_{\delta'}[A]$,
	\item for each $x \in L_{\delta'}[A]$ there is a bijection $f \in L_{\delta'}[A]$ between $\omega_1$ and $x$,
	\item \label{d'3} taking the sentence $\sigma$ from Lemma $\ref{conL}$ $(L_{\delta'}[A], \in,  A) \models \sigma$.
\end{enumerate} 
\ed
(Using Claim $\ref{haszncor}$ and $(|L_\alpha[A]| = |\alpha|)^{L[A]}$ for $\alpha \geq \omega$ it is easy to see that we can do this closure operation, and there is such a $\delta'<\omega_2$.) 
Then we have
\beeq \label{doo}\left( \delta' \text{ is a limit } \right)  \bigwedge \ \left( L_{\delta'}[A] \models \ \text{ '} \omega_1 \text{ is the largest cardinal'} \right), \eeq
and also the desired uniqueness by our next claim.
\bc \label{sigma''}There is a statement $\sigma' \in \mathcal{L}_\in(R_A, c_{\partial})$ such that for each $\delta \in [\xi, \omega_2)$ $(L_{\delta'}[A], \in, A, \delta) \models \sigma'$, moreover, for each $ \delta > \omega_1$ and $\delta'' > \delta$ 
\[ ((L_{\delta''}[A], \in, A, \delta) \models \sigma') \ \Rightarrow \  (\delta'' = \delta').\]
\ec
\bp
 First define $\sigma'' = \sigma \wedge (\forall y \exists f: \omega_1 \to y \text{ bijection})$ and let $\sigma'$ be the following sentence
 \[ \sigma' = \sigma'' \wedge \left(\neg(\exists X) \ (X \text{ is transitive}) \wedge (\sigma'')^X \wedge (\delta \in X) \right)  \]
 (where under $\psi^X$ we always mean the formula $\psi \in \mathcal{L}_\in(R_A, c_\partial)$ relativized to $X$, and $\sigma$ is from Lemma $\ref{conL}$).
\ep

Now fix $\delta \in [\xi, \omega_2)$, and for each ordinal $0<\alpha < \omega_1$ define 
$M_{\delta, \alpha}$ to be the Skolem-hull  
\beeq \label{mdf} M_{\delta, \alpha} =\mathfrak{H}^{(L_{\delta'}[A], \in, A, \delta)}(\alpha) \ \ (\text{for each  }\alpha < \omega_1), \eeq 
Also define 
\beeq M_{\delta, 0} = \emptyset. \eeq
Then 
\beeq \label{elemi} M_{\delta, \alpha} \prec ( L_{\delta'}[A], \in , A, \delta ) \ \ (\text{for each } \alpha>0). \eeq  
Observe that whenever $M^* \prec (L_{\delta'}[A], \in, A, \delta)$  we have for the Skolem functions from Definition $\ref{Skdef}$ that $f^{(L_{\delta'}[A], \in, A, \delta)}_\psi \up (M^*)^{n_\psi} = f^{(M^*, \in, A \cap M^*, \delta)}_\psi$ ,
hence
\beeq \label{Skabsz} \forall M' \subseteq M^* \prec (L_{\delta'}[A], \in, A, \delta): \ \ \mathfrak{H}^{(L_{\delta'}[A], \in, A, \delta)}(M') = \mathfrak{H}^{(M^*, \in, A \cap M^*, \delta)}(M'). \eeq


Now as we defined $\langle M_{\delta, \alpha}: \ \alpha < \omega_1 \rangle$ note that \beeq \label{kezdosze} (M \prec (L_{\delta'}[A], \in, A, \delta)) \wedge (|M| = \omega)  \ \to \ (M \cap \omega_1 \in \omega_1), \eeq
in particular
\beeq  M_{\delta, \alpha} \cap \omega_1 \in \omega_1, \eeq
 since $\eqref{xi}$ together with $\xi \leq \delta < \delta'$ implies that in $L_{\delta'}[A]$ there is an enumeration of each ordinal less than $\omega_1$ (and $M_{\delta, \alpha}$ is countable).
This implies that
\[ (C_{\delta} = \{ \alpha < \omega_1: \ M_{\delta,\alpha} \cap \omega_1 = \alpha \} \text{ is a club in } \omega_1) \ \wedge \ (0 \in C_{\delta}). \]
It is easy to see that 
\beeq \label{Cd} \forall \alpha < \omega_1: \ M_{\delta,\alpha} = M_{\delta, \min(C_\delta \sm \alpha)}. \eeq

For later use we verify the following statement.
\bc \label{kiad}
 \[ \bigcup_{\alpha < \omega_1} M_{\delta, \alpha} = L_{\delta'}[A]. \]
\ec
\bp
Since the union of an increasing chain of elementary submodels	is an elementary submodel, we have $M_{\omega_1} = \bigcup_{\alpha < \omega_1} M_{\delta, \alpha} \prec (L_{\delta'}[A], \in, A, \delta)$. Now recall, that in $L_{\delta'}[A]$ every set $x$ admits a surjection from $\omega_1$ onto $x$, therefore $\omega_1 \subseteq M_{\omega_1}$ implies that $M_{\omega_1}$ is transitive.
Then by Lemma $\ref{conL}$ and $M_{\omega_1} \models \sigma$ we have  $M_{\omega_1} = L_{\delta''}[A]$ for some $\delta'' > \delta$. But then either $M_{\omega_1} \in L_{\delta'}[A]$, or $M_{\omega_1} = L_{\delta'}[A]$, and because the former would contradict Claim $\ref{sigma''}$, we arrive at our conclusion.
\ep

For each $\alpha \in C_{\delta}$ and $\beta< \omega_1$, if $\alpha = \max (C_{\delta} \cap (\beta+1))$, then let $N_{\delta, \beta, \alpha}$ be the range of the Mostowski-collapse $\pi_{\delta, \alpha}$ of $(M_{\delta, \alpha}, \in)$,
and let $A_{\delta, \beta, \alpha} = \pi_{\delta, \alpha}(A)$, $\partial_{\delta, \beta, \alpha} = \pi_{\delta, \alpha}(\delta)$:
\beeq \pi_{\delta, \alpha}: M_{\delta, \alpha} \to N_{\delta, \beta, \alpha}, \eeq
which is of course not only an isomorphism between $(M_{\delta, \alpha}, \in)$ and  $(N_{\delta, \beta, \alpha}, \in)$, but witnesses
\beeq \label{extiso} (M_{\delta, \alpha}, \in, A  \cap M_{\delta, \alpha}, \delta )  \simeq (N_{\delta, \beta, \alpha}, \in , A_{\delta, \beta, \alpha}, \partial_{\delta, \beta, \alpha}). \eeq




Now we are ready to construct the tree $T$.
For a fixed $\delta \in [\xi, \omega_2)$, $\alpha \in C_{\delta}$, $\beta < \omega_1$, if $0< \alpha = \max(C_{\delta} \cap (\beta +1))$ holds then we define
\beeq t_{\delta, \beta, \alpha} = (N_{\delta, \beta, \alpha}, \in, A_{\delta, \beta, \alpha}, \partial_{\delta, \beta, \alpha}), \eeq
i.e. the structure $(N_{\delta, \beta, \alpha}, \in)$ extended by the one-place relation for the image of $A \in M_{\delta, \alpha}$ under the collapsing isomorphism, and the constant symbol for $\partial_{\delta, \beta, \alpha}$. For $\max(C_{\delta} \cap (\beta +1))=0$ let $ t_{\delta, \beta, 0} = \emptyset$.
\\ Observe that given $t = t_{\delta, \beta, \alpha}$ we can decode $\alpha$ from $t$, as $\alpha$ is the first uncountable ordinal of $t$.

\bd \label{rend}
Define
\[ T = \{ (\beta, t_{\delta,\beta, \alpha}): \ \delta \in [\xi,\omega_2), \beta < \omega_1, \ \alpha = \max(C_{\delta} \cap (\beta+1)) \}, \]
with the partial order $(\beta_0, t_{\delta_0,\beta_0, \alpha_0}) \leq_{T}  (\beta_1, t_{\delta_1, \beta_1, \alpha_1})$ iff  either $\alpha_0 = 0$ (thus $t_{\delta_0, \beta_0, \alpha_0}$ is the empty structure), or
\begin{enumerate}[label = (\roman*)]
	\item $ \beta_0 \leq \beta_1$, and 
	 \item \label{ket} taking the Skolem-hull $M$ of $\alpha_0$ in 
	 \[t_{\delta_1,\beta_1, \alpha_1} = (N_{\delta_1, \beta_1, \alpha_1}, \in, A_{\delta_1, \beta_1, \alpha_1}, \partial_{\delta_1, \beta_1, \alpha_1}) \] (i.e. $M = \mathfrak{H}^{t_{\delta_1,\beta_1, \alpha_1}}(\alpha_0)$ is isomorphic to 
	 $t_{\delta_0,\beta_0, \alpha_0}$:
	 \[  (M, \in, A_{\delta_1,\beta_1, \alpha_1} \cap M, \partial_{\delta_1,\beta_1, \alpha_1}) \simeq ( N_{\delta_0, \beta_0, \alpha_0}, \in, A_{\delta_0,\beta_0, \alpha_0}, \partial_{\delta_0,\beta_0, \alpha_0}),\] 
	 and
	 
	 \item \label{ha} if $\alpha_0 < \alpha_1$, then there is no proper elementary submodel $M \prec (N_{\delta_1, \beta_1, \alpha_1}, \in, A_{\delta_1,\beta_1, \alpha_1}, \partial_{\delta_1,\beta_1, \alpha_1})$ with
	 \[ \alpha_0 \cup \{\alpha_0\} \subseteq M, \text{ and} \]
	 \[  M \cap \alpha_1 \subseteq \beta_0.\] 
\end{enumerate}
\ed
Roughly speaking, in level $\beta$ we have (isomorphism types of) initial segments $M$ of models of the form $(L_{\Delta'}[A], \in, A, \Delta)$ (for some $\Delta \in [\xi,\omega_2)$), such that $M \cap \omega_1 \leq \beta$, and $M$ is maximal w.r.t. this condition.
We need to check that $T$ is a tree, its levels are countable, and that it has only $\omega_2$-many branches even in $V$.

The following claim is a standard calculation, but for the sake of completeness we include the proof.
\bc \label{uni} Let $\delta \in [\xi, \omega_2)$ be fixed, $\beta_0 \leq \beta_1 < \omega_1$,
let $\alpha_1 = \max(C_{\delta} \cap (\beta_1+1))$, $\alpha_0 = \max(C_{\delta} \cap (\beta_0+1))$. Then $(\beta_0, t_{\delta,\beta_0, \alpha_0}) \leq_{T}  (\beta_1, t_{\delta, \beta_1, \alpha_1})$.

Moreover, the embedding $\pi_{\beta_0, \beta_1}: N_{\delta,\beta_0, \alpha_0} \to N_{\delta,\beta_1, \alpha_1}$ is unique.
\ec
\bp
First observe that by $\eqref{mdf}$ and $\eqref{Skabsz}$ for $\delta \in [\xi, \omega_2)$, $\alpha_0 < \alpha_1$
\[ \mathfrak{H}^{(M_{\delta, \alpha_1},\in, A, \delta)} ( \alpha_0) = \mathfrak{H}^{(L_{\delta'}[A], \in, A, \delta)} ( \alpha_0) = M_{\delta, \alpha_0}, \]
therefore since $\beta_1< \omega_1$ is such that $\alpha_1 = \max( C_{\delta} \cap (\beta_1+1))$, then applying (the restriction of) the collapsing isomorphism $\pi_{\delta, \alpha_1}$ to the left side, we obtain
\[ (\mathfrak{H}^{(N_{\delta,  \beta_1,\alpha_1}, \in, A_{\delta, \beta_1, \alpha_1}, \partial_{\delta,\beta_1, \alpha_1})} ( \alpha_0 ), \in) \simeq  (M_{\delta, \alpha_0}, \in) \]
and because $\beta_0 < \beta_1$ is such that $\alpha_0 = \max( C_{\delta} \cap (\beta_0+1))$, then applying the isomorphism $\pi_{\delta, \alpha_0}$ to the right side (which fixes $\alpha_0$) we obtain
\[ (\mathfrak{H}^{(N_{\delta,  \beta_1,\alpha_1}, \in, A_{\delta, \beta_1, \alpha_1}, \partial_{\delta,\beta_1, \alpha_1})} ( \alpha_0 ), \in) \simeq (N_{\delta, \alpha_0, \beta_0}, \in). \]
Finally, since $\pi_{\delta, \alpha_1}(A) = A_{\delta,\beta_1, \alpha_1}$, $\pi_{\delta, \alpha_0}(A) = A_{\delta,\beta_0, \alpha_0}$, and
$\pi_{\delta, \alpha_1}(\delta) = \partial_{\delta,\beta_1, \alpha_1}$, $\pi_{\delta, \alpha_0}(\delta) = \partial_{\delta,\beta_0, \alpha_0}$, we have
\[ (\mathfrak{H}^{N_{\delta, \beta_1, \alpha_1}} ( \alpha_0), \in  A_{\delta,\beta_1, \alpha_1}, \partial_{\delta,\beta_1, \alpha_1})  \]
\[ \text{is isomorphic to } (N_{\delta, \beta_0, \alpha_0}, \in,  A_{\delta,\beta_0, \alpha_0}, \partial_{\delta,\beta_0, \alpha_0}), \]
therefore $\ref{ket}$ holds. The uniqueness easily follows from the facts that the embedding of $(N_{\delta, \beta_0, \alpha_0}, \in,  A_{\delta,\beta_0, \alpha_0}, \partial_{\delta,\beta_0, \alpha_0})$ has to fix the ordinals less than $\alpha_0$, and elementary embeddings uniquely extend to Skolem-hulls.

For $\ref{ha}$ suppose that $\alpha_0 < \alpha_1$, and note that
\[ (N_{\delta, \beta_1, \alpha_1}, \in) \models \ \text{'} \alpha_1 \text{ is the least uncountable ordinal, }  \alpha_0 \text{ is countable',} \]
and for $M \prec (N_{\delta, \beta_1, \alpha_1}, \in, A_{\delta, \beta_1, \alpha_1} ,\partial_{\delta, \beta_1, \alpha_1})$ if $\alpha_0 \cup \{\alpha_0\} \subseteq M$ then consider the corresponding submodel $M' \prec (M_{\delta,\alpha_1}, \in, A, \delta)$, for which $M' \supseteq M_{\delta, \alpha_0+1}$. But (recalling $\eqref{kezdosze}$) since $\max(C_{\delta} \cap (\beta_0+1)) = \alpha_0$ we obtain $\beta_0 \cup \{\beta_0 \} \subseteq M' \subseteq M_{\delta, \alpha_1}$, that can happen only if $\beta_0$ is smaller than the least uncountable ordinal in $N_{\delta, \beta_1, \alpha_1}$,  $\alpha_1$. But then $\beta_0 \in M \cap \alpha_1$.
\ep

The next claim will verify that $T$ is a tree of height $\omega_1$ (for the transitivity of $\leq_T$ use the claim two times).

\bc \label{acsa} For a fixed $\delta_1 \in [\xi, \omega_2)$, $\beta_0 \leq \beta_1 < \omega_1$,
let $\alpha_1 = \max(C_{\delta_1} \cap (\beta_1+1)$, and fix arbitrary $\alpha_0 \in \omega_1$, $\delta_0 \in [\xi, \omega_2)$. Then $(\beta_0, t_{\delta_0,\beta_0, \alpha_0}) \leq_{T}  (\beta_1,t_{\delta_1, \beta_1, \alpha_1})$ iff $t_{\delta_0,\beta_0, \alpha_0} = t_{\delta_1, \beta_0, \max(C_{\delta_1 \cap (\beta_0+1) })}$.
\ec
\bp
We only have to check the 'only if' part, but first observe that Definition $\ref{rend}$ clearly implies that up to isomorphism there exists only one $t$ for which $(\beta_0, t) \leq (\beta_1,t_{\delta_1, \beta_1, \alpha_1})$.  Now the claim is the consequence of the fact that 
$t_{\delta_*,\beta_0, \alpha_*} \neq t_{\delta_{**}, \beta_0, \alpha_{**}}$ implies that they are not isomorphic as  structures of the language $\mathcal{L}_\in(R_A, c_\partial)$:
 For transitive sets $N$ and $N'$ with $X, \partial \in N$, $X', \partial' \in N'$ the structures $(N, \in, X, \partial)$, $(N', \in, X', \partial')$ are isomorphic if and only if $N = N'$, $X = X'$ and $\partial = \partial'$ (since by the uniqueness of the Mostowski collapse we know that $(N, \in) \simeq (N', \in)$ iff $N=N'$).
\ep

\bl \label{blev} For each $\beta < \omega_1$ the $\beta$'th level of $T$ is countable.
\el
\bp
By Claim $\ref{acsa}$ we have that the $\beta$'th level of $T$ is
\[ T_{\leq \beta} \sm T_{<\beta} =\{ (\beta,t_{\delta,\beta, \alpha}): \ \delta \in [\xi,\omega_2), \ \alpha = \max(C_{\delta} \cap (\beta+1)) \} \}.\]
For a fixed $\delta \in [\xi,\omega_2)$ fix $\alpha = \max(C_{\delta} \cap (\beta+1))$ too, and consider the structure 
\[ t_{\delta,\beta, \alpha} = (N_{\delta, \beta, \alpha}, \in, A_{\delta, \beta, \alpha} , \partial_{\delta, \beta, \alpha} ), \]
where $N_{\delta, \beta, \alpha}$ is the Mostowski collapse of $(M_{\delta, \alpha}, \in)$ (by the isomorphism $\pi_{\delta,\alpha}$), and $A_{\delta, \beta, \alpha} = A \cap \alpha$. Now 
$\eqref{elemi}$ states $M_{\delta, \alpha} \prec (L_{\delta'}, \in , A)$
then (recalling $M_{\delta, \alpha} \cap \omega_1 = \alpha$, and $\pi_{\delta, \alpha} \up \alpha = \id_{\alpha}$) by Lemma $\ref{conL}$
\[ N_{\delta, \beta, \alpha} = L_{\gamma}[A \cap \alpha] \]
for some $\gamma = \gamma(\delta, \alpha) \in (\alpha,\omega_1)$. 
Now we determine an upper bound $\gamma_\alpha$ for the set $\{\gamma(\delta, \alpha): \delta \in [\xi,\omega_2) \ \wedge \  \alpha \in C_{\delta} \}$. If we have such a bound for each possible $\alpha \leq \beta$, then letting $\gamma_\infty$ denote $\sup\{\gamma_\alpha: \ \alpha \leq \beta\}$, we get 
\[ \{ t_{\delta,\beta, \alpha}): \ \delta \in [\xi,\omega_2), \ \alpha = \max(C_{\delta} \cap (\beta+1)) \} \} \subseteq \]
\[ \{ (L_\gamma[A\cap \alpha], \in, A \cap \alpha,  \partial ): \ \gamma \leq \gamma_\infty, \ \alpha \leq \beta, \ \partial < \gamma \}, \]
which latter set is obviously countable, this will finish the proof of the lemma.

So fix $\alpha \leq \beta$ and $\delta \in [\xi, \omega_2)$ such that $\alpha \in C_{\delta}$. Now we have two cases depending on whether there is any (cardinal)$^{L[A \cap \alpha]}$ in $(\alpha, \omega_1)$.
If $\lambda \in (\alpha, \omega_1)$ is a cardinal in the inner model $L[A \cap \alpha]$, then for each $\delta$ if $\alpha = \max(C_{\delta} \cap (\beta+1))$, then the transitive set $N_{\delta, \beta, \alpha}$ cannot contain $\lambda$, as  $M_{\delta, \alpha}$ sees $\omega_1$ as the largest cardinal, and $\pi_{\delta, \alpha}(\omega_1) = \alpha$. This case choosing $\gamma_\alpha = \lambda$ works.

On the other hand, if $(|\alpha|^+)^{L[A \cap \alpha]} = \omega_1$, then we first prove that $\alpha \in C_{\delta}$ implies $(|\alpha| = \omega)^{L[A \cap \alpha]}$: 
otherwise in $M_{\delta, \alpha}$, and in $N_{\delta, \beta, \alpha}$ each ordinal less than $\alpha$ are countable, thus as well in $L[A \cap \alpha]$. Then it is easy to see that the condition
\[ (\lambda \text{ is the unique cardinal in } (\omega, \omega_1^V))^{L[A \cap \lambda]}\]
cannot hold for two different $\lambda$'s, therefore $\alpha$ can be defined in $L[A]$. But then using Claim $\ref{haszncor}$ with $X = A \cap \alpha$ we have that for each $\zeta \in (\alpha, \omega_1)$ there is a bijection $f_\zeta \in L_{\omega_1}[A \cap \alpha]$ between $\alpha$ and $\zeta$, therefore $\alpha$ can be defined also in $L_{\delta'}[A]$, and $M \prec (L_{\delta'}[A], \in)$ implies $\alpha \in M$, contradicting that $M_{\delta, \alpha} \cap \omega_1 = \alpha$ (which holds by $\alpha \in C_{\delta}$).
Then $(|\alpha| = \omega)^{L[A \cap \alpha]}$ and Claim $\ref{haszncor}$ implies that there is an ordinal $\lambda < \omega_1$ such that there exists a bijection between $\alpha$ and $\omega$ in $L_{\lambda}[A \cap \alpha]$, implying
\[ N_{\delta, \beta, \alpha} = L_{\gamma(\delta, \alpha)}[A \cap \alpha] \subsetneq L_{\lambda}[A \cap \alpha], \]
since $\alpha$ is uncountable in $N_{\delta, \beta, \alpha}$. This case 
\[ \{\gamma(\delta, \alpha): \delta \in [\xi,\omega_2) \ \wedge \  \alpha \in C_{\delta} \} \subseteq \gamma_\alpha = \lambda, \]
which completes the proof of Lemma $\ref{blev}$.

\ep

Now $T$ is obviously a Kurepa tree by the following fact and lemma.
\bfa The sequence $\langle B_\delta: \ \delta \in [\xi, \omega_2) \rangle$ lists pairwise distinct cofinal branches in $T$, where 
\[ B_{\delta}= \{ (\beta,t_{\delta, \beta, \max(C_{\delta} \cap (\beta+1))}): \ \beta < \omega_1 \}.  \]

\efa
\bp
We only need to prove that $B_\delta \neq B_\gamma$ if $\delta \neq \gamma$. But according to the second statement of Claim $\ref{uni}$ for each $\beta < \beta' < \omega_1$ there is a unique elementary embedding of $t_{\delta,  \beta', \max(C_{\delta} \cap (\beta'+1))}$ to  $t_{\delta,  \beta, \max(C_{\delta} \cap (\beta+1))}$, therefore there is a unique direct-limit of this elementary chain, isomorphic to $\bigcup_{\alpha \in C_\delta} M_{\delta, \alpha}$, which is $(L_{\delta'}[A], \in, A, \delta)$ by Claim $\ref{kiad}$. 
\ep

It is only left to prove that  each branch of $T$ is of the form  $B_{\delta}$ for some $\delta \in [\xi, \omega_2)$ (even in $V$). The following lemma will complete the proof of Theorem $\ref{szukst}$.

\bl
Let $B \subseteq T$ a cofinal branch in $T$, $B \in V$.
Then $B = B_{\delta_\bullet}$ for a unique $\delta_\bullet \in [\xi, \omega_2)$.
\el
\bp
Let $t_{\delta_\beta, \beta, \alpha_\beta} = (N_{\delta_\beta, \beta, \alpha_\beta}, \in, A_{\delta_\beta, \beta, \alpha_\beta}, \partial_{\delta_\beta, \beta, \alpha_\beta})$ denote the element in $B \cap (T_{\leq \beta} \sm T_{< \beta})$.
Working in $V$ first we define the following bonding maps: for $\gamma \leq \beta < \omega_1$ let 
\[ \pi_{\gamma,\beta}: N_{\delta_\gamma, \gamma, \alpha_\gamma} \to N_{\delta_\beta, \beta, \alpha_\beta} \] 
be the unique elementary embedding (combining Claim $\ref{acsa}$, and the second statement of Claim $\ref{uni}$). Since elementary submodels of an elementary submodel are elementary submodels, $\pi_{\beta', \beta} \circ \pi_{\beta'', \beta'}$ is an elementary embedding for each $\beta'' \leq \beta' \leq \beta < \omega_1$, therefore by the uniqueness 
\beeq \label{unik} (\forall \beta'' \leq \beta' \leq \beta < \omega_1): \  \pi_{\beta', \beta} \circ \pi_{\beta'', \beta'} = \pi_{\beta'', \beta}. \eeq
This elementary chain allows us to define the limit $D = (N_{\omega_1},  \mathbf{E}, A_{\omega_1}, \partial_{\omega_1})$ of the directed system $\{t_{\delta_\beta, \beta, \alpha_\beta}, \ \pi_{\beta',\beta}: \ \ \beta' \leq \beta < \omega_1 \}$. 

Let $\pi_\beta: N_{\delta_\beta, \beta, \alpha_\beta} \to N_{\omega_1}$ be the embedding, $N_\beta = \ran(\pi_\beta)$ (hence $N_{\omega_1} = \bigcup_{\beta< \omega_1} N_\beta$).

First note that $(N_{\omega_1}, \mathbf{E})$ is well-founded, otherwise there would be an infinite $\mathbf{E}$-decreasing  chain in the embedded image of $N_{\delta_\beta, \beta, \alpha_\beta}$ for some (in fact, every large enough) $\beta$, contradicting that $(N_{\delta_\beta, \beta, \alpha_\beta}, \in)$ is well-founded.
Now (by the $\mathbf{E}$-extensionality in $N_{\omega_1}$) we can assume that $N_{\omega_1}$ is a Mostowski collapse, i.e. $ (N_{\omega_1} , \mathbf{E}) = (N_{\omega_1}, \in)$. Then it is easy to see that if $\beta< \omega_1$ for the elementary embedding $\pi_\beta: N_{\delta_\beta, \beta, \alpha_\beta} \to N_{\omega_1}$ we have $\pi_\beta \up \alpha_\beta = \id_{\alpha_\beta}$, and $\pi_\beta(\alpha_\beta) = \omega_1$, thus (recalling that $A_{\delta_\beta, \beta, \alpha_\beta} = A \cap \alpha_\beta$) we obtain $(N_{\omega_1},  \mathbf{E}, A_{\omega_1}, \partial_{\omega_1}) = (N_{\omega_1}, \in, A, \delta_{\bullet})$ for some $\delta_\bullet \in (\omega_1, \omega_2)$.
Now we can use Lemma $\ref{conL}$ (since $(N_{\delta_\beta, \beta, \alpha_\beta}, \in, A_{\delta_\beta, \beta, \alpha_\beta}) \models \sigma$), there exists $\zeta> \delta_\bullet$ such that
\[ N_{\omega_1} = L_{\zeta}[A], \]
and then
\[ (N_{\omega_1}, \in, A, \delta_\bullet) = (L_{\zeta}[A], \in, A, \delta_\bullet). \]
Now because the formula $\sigma' \in \mathcal{L}_{\in}(R_A, c_\partial)$ from Claim $\ref{sigma''}$ holds in $(L_{\delta'}[A], \in, A, \delta)$ (for each $\delta \in [\xi, \omega_2)$) (for our mapping $\delta \mapsto \delta'$ from Definition $\ref{delta'}$) and therefore also in $M_{\delta, \alpha}$'s, $N_{\delta, \beta, \alpha}$'s ($\delta \in [\xi, \omega_2)$), so it must hold in $(N_{\omega_1}, \in, A, \delta_\bullet)$, which means that $\delta_\bullet \geq \xi$, and $\zeta = \delta_\bullet'$ , i.e.
\[ (N_{\omega_1}, \in, A, \delta_\bullet) = (L_{\delta_\bullet'}[A], \in, A, \delta_\bullet), \]
Finally, we have to prove that for each $\beta < \omega_1$ 
\[ t_{\delta_\beta, \beta, \alpha_\beta} = (N_{\delta_\beta, \beta, \alpha_\beta}, \in, A_{\delta_\beta, \beta, \alpha_\beta}, \partial_{\delta_\beta, \beta, \alpha_\beta}) = t_{\delta_\bullet, \beta, \max(C_{\delta_\bullet} \cap (\beta+1))} \]
by arguing (having $\beta$ fixed) that for a large enough $\gamma$  
\[ (\beta,t_{\delta_\bullet, \beta, \max(C_{\delta_\bullet} \cap (\beta+1))}) \leq_T (\gamma, t_{\delta_\gamma, \gamma, \alpha_\gamma}). \]
Let $\alpha = \max(C_{\delta_\bullet} \cap (\beta+1))$, $\alpha' = \min (C_{\delta_\bullet} \sm (\beta+1))$, $\beta' = \alpha'$, and consider the models $M_{\delta_\bullet, \alpha}, M_{\delta_\bullet, \alpha'} \prec (L_{\delta_\bullet'}[A], \in, A, \delta_\bullet)$. Choose $\gamma \geq \beta'$, $\gamma < \omega_1$ so that $N_{\gamma
	} = \pi_\gamma[N_{\delta_\gamma, \gamma, \alpha_\gamma}]  \supseteq M_{\delta_\bullet, \alpha'}$. Then
	\beeq \label{alf} \alpha_\gamma \geq \alpha' > \beta +1, \text{ and} \eeq
	 $\alpha' \cup \{ \omega_1 \} \subseteq N_{\gamma} \prec (L_{\delta_\bullet'}[A], \in, A, \delta_\bullet)$  with $\eqref{Skabsz}$ imply
	\[ \mathfrak{H}^{(N_{\gamma}, \in, A \cap N_{\gamma}, \delta_\bullet)}(\alpha) = \mathfrak{H}^{(L_{\delta_\bullet'}[A], \in, A, \delta_\bullet)} (\alpha) = M_{\delta_\bullet, \alpha}. \]
	Therefore in $ (N_{\gamma}, \in, A \cap N_{\gamma}, \delta_\bullet) \simeq (N_{\delta_\gamma, \gamma, \alpha_\gamma}, \in, A_{\delta_\gamma, \gamma, \alpha_\gamma}, \partial_{\delta_\gamma, \gamma, \alpha_\gamma})$ there is an elementary submodel isomorphic to $(M_{\delta_\bullet, \alpha}, \in, A \cap M_{{\delta_\bullet}, \alpha}, \delta_\bullet)$, which latter is isomorphic to $(N_{\delta_\bullet,  \beta, \alpha}, \in, A \cap \alpha, \partial_{\delta_\bullet, \beta, \alpha})$, thus $\ref{ket}$ from Definition $\ref{rend}$ holds.
 \\ Similarly, using also $\eqref{Cd}$ and the definitions of $\alpha$, $\alpha'$
\[ \mathfrak{H}^{(N_{\gamma
	}, \in, A \cap N_{\gamma}, \delta_\bullet)}(\alpha+1) = M_{\delta_\bullet, \alpha+1} = M_{\delta_\bullet, \alpha'} \supseteq \alpha' \supseteq \beta \cup \{ \beta\},\]
ans since the isomorphism between $(N_{\gamma
}, \in, A \cap N_{\gamma
}, \delta_\bullet)$ and $(N_{\delta_\gamma, \gamma, \alpha_\gamma}, \in, A_{\delta_\gamma, \gamma, \alpha_\gamma}, \partial_{\delta_\gamma, \gamma, \alpha_\gamma})$ fixes the ordinals less than or equal to $\alpha'$ we obtain
\[  \mathfrak{H}^{(N_{\delta_\gamma, \gamma, \alpha_\gamma}, \in, A_{\delta_\gamma, \gamma, \alpha_\gamma}, \partial_{\delta_\gamma, \gamma, \alpha_\gamma}})(\alpha+1) \supseteq \beta \cup \{ \beta\}. \] 
Therefore recalling $\eqref{alf}$ we obtain that $\ref{ha}$ (of Definition $\ref{rend}$) holds as well. 

\ep



\bibliographystyle{amsalpha}
\bibliography{listb,listx,1189}

%
%
%
%
%
%
%
%

\end{document}